\documentclass[english]{article}
\usepackage[T1]{fontenc}
\usepackage[latin9]{inputenc}
\usepackage{geometry}
\usepackage{xcolor}
\usepackage{babel}
\usepackage{verbatim}
\usepackage{prettyref}
\usepackage{float}
\usepackage{mathtools}
\usepackage{algorithm2e}
\usepackage{amsmath}
\usepackage{amsthm}
\usepackage{amssymb}
\usepackage{graphicx}
\usepackage{subfig}
\usepackage[authoryear]{natbib}
\usepackage[unicode=true,pdfusetitle,
 bookmarks=true,bookmarksnumbered=false,bookmarksopen=false,
 breaklinks=true,pdfborder={0 0 1},backref=false,colorlinks=true]
 {hyperref}
\hypersetup{
 linkcolor=blue,citecolor=blue,urlcolor=blue,filecolor=blue,pdfpagelayout=OneColumn,pdfnewwindow=true,pdfstartview=XYZ,plainpages=false,pdfpagelabels,hyperindex=true}

\makeatletter

\floatstyle{ruled}
\newfloat{algorithm}{tbp}{loa}
\providecommand{\algorithmname}{Algorithm}
\floatname{algorithm}{\protect\algorithmname}

\theoremstyle{definition}
    \ifx\thechapter\undefined
      \newtheorem{defn}{\protect\definitionname}
    \else
      \newtheorem{defn}{\protect\definitionname}[chapter]
    \fi
\theoremstyle{plain}
    \ifx\thechapter\undefined
      \newtheorem{assumption}{\protect\assumptionname}
    \else
      \newtheorem{assumption}{\protect\assumptionname}[chapter]
    \fi
\theoremstyle{plain}
    \ifx\thechapter\undefined
      \newtheorem{prop}{\protect\propositionname}
    \else
      
    \fi
\theoremstyle{plain}
    \ifx\thechapter\undefined
	    \newtheorem{thm}{\protect\theoremname}
	  \else
      \newtheorem{thm}{\protect\theoremname}[chapter]
    \fi
\theoremstyle{plain}
    \ifx\thechapter\undefined
  \newtheorem{cor}{\protect\corollaryname}
\else
      
    \fi
\theoremstyle{plain}
    \ifx\thechapter\undefined
  \newtheorem{prob}{\protect\problemname}
\else
      
    \fi
\theoremstyle{plain}
    \ifx\thechapter\undefined
  \newtheorem{remark}{\protect\remarkname}
\else
      
    \fi

\@ifundefined{date}{}{\date{}}
\usepackage{babel}

\usepackage{stackengine}
\usepackage{bbm}
\usepackage{enumitem}

\usepackage{authblk}

\setlist{leftmargin=*, topsep=0.5em, parsep=0pt, itemsep=1em, labelindent=0pt, align=left}

\@ifundefined{showcaptionsetup}{}{%
 \PassOptionsToPackage{caption=false}{subfig}}
\usepackage{subfig}
\makeatother

\providecommand{\assumptionname}{Assumption}
\providecommand{\corollaryname}{Corollary}
\providecommand{\definitionname}{Definition}
\providecommand{\propositionname}{Proposition}
\providecommand{\remarkname}{Remark}
\providecommand{\theoremname}{Theorem}
\providecommand{\problemname}{Problem}

\usepackage{nameref}

\newtheorem{problem}{Problem}

\begin{document}

\title{Deep Semi-Martingale Optimal Transport}
\author[1,2]{Ivan Guo}
\author[3]{Nicolas Langren\'e}
\author[1,2,4]{Gr\'egoire Loeper}
\author[1]{Wei Ning}

\affil[1]{\small School of Mathematical Sciences, Monash University, Melbourne, Australia}
\affil[2]{\small Centre for Quantitative Finance and Investment Strategies, Monash University, Australia}
\affil[3]{\small Data61, Commonwealth Scientific and Industrial Research Organisation, Australia}
\affil[4]{\small BNP Paribas Global Markets}

\maketitle

\begin{abstract}
We propose two deep neural network-based methods for solving semi-martingale optimal transport problems. 
The first method is based on a relaxation/penalization of the terminal constraint, and is solved using deep neural networks.
The second method is based on the dual formulation of the problem, which we express as a saddle point problem,
and is solved using adversarial networks. Both methods are mesh-free and therefore mitigate the {\it curse of dimensionality}.
We test the performance and accuracy of our methods on several examples up
to dimension 10.
We also apply the first algorithm to a portfolio optimization problem where the goal is, given an initial wealth distribution, to find an investment strategy leading to a prescribed terminal wealth distribution.
\end{abstract}

\section{Introduction}

The optimal transport problem goes back to
the work of \citet{monge1781memoire} and aims at transporting a distribution 
$\mu$ to another distribution $\nu$ under minimum transport cost.
It was later revisited by \citet{kantorovich1942translocation},  leading to the so-called Monge--Kantorovich formulation. 
In recent years, a fast-developing phase was spurred by a wide range of extensions and applications of the Monge--Kantorovich problem; interested readers can refer to the books by \citet{rachev1998mass}, \citet{villani2003topics} and \citet{villani2008optimal} for a comprehensive review.
Although we have gained tremendous theoretical insight, the numerical
solution of the problem remains challenging. When the dimension is less or equal
to three, many state-of-art approaches are able
to compute the global solution effectively; see, for example, \citet{chow2019algorithm}, \citet{haber2015multilevel}, \citet{li2018parallel}, and the review by \citet{zhang2020review}. Readers can refer to the books by \citet{santambrogio2015optimal} and \citet{peyre2019computational}, and the references therein for an overview of these approaches.
However, many traditional methods rely on Euclidean coordinates and require spatial discretization. When the distributions live in spaces of dimension four or more, these traditional
methods suffer from the curse of dimensionality. Under this situation,
solving optimal transport problems using deep neural networks looks very
attractive since it can avoid space discretization.

Artificial Neural Networks (ANNs) are at the core of the Deep Learning methodology. ANNs were first introduced long ago by \citet{mcculloch1943logical},
but due to
the tremendous increase in computing power they have been widely used in the recent years. When an ANN has two or
more hidden layers, it is called a Deep Neural Network (DNN). \footnote{DNNs are preferred
to shallow networks in general. As shown in \citet{liang2016deep},
for a given degree of approximation error, the number of neurons needed
by a shallow network to approximate the function is exponentially
larger than the number of neurons needed by a deep network. }
Usually the calibration of an ANN is done by minimizing a loss function.  When the loss function is an expectation,
 Stochastic Gradient Descent (SGD) is a natural adaptation of Gradient Descent to stochastic
optimization problems. 
The DNN-SGD paradigm is well suited for solving scientific computing problems such as 
stochastic control problems (e.g., \citealt{han2016deep}, \citealt{hure2018deep}, \citealt{bachouch2018deep})
or partial differential equations (PDEs) in very high dimension (e.g., \citealt{weinan2017deep},
\citealt{sirignano2018dgm} and \citealt{hure2019some}), thanks to its ability to overcome the curse of dimensionality. 
A comprehensive review of the numerical and theoretical advances in solving PDE and BSDE with deep learning
algorithms can be found in \citet{han2020algorithms}.

Recent years witnessed the emergence of research on solving optimal transport
problems with neural networks. \citet{ruthotto2020machine} used a
residual network (ResNet) to approximately solve high-dimensional
mean-field games by combining Lagrangian and Eulerian viewpoints. In
the numerical experiment, they solved a dynamical optimal transport problem
as a potential mean-field game. \citet{henry2019martingale} introduced
Lagrange multipliers associated with the two marginal constraints and
proposed a Lagrangian algorithm to solve martingale optimal transport
using neural networks. However, for the sake of simplicity, this work only focuses on cost functions satisfying a martingale condition. \citet{eckstein2019computation}
studied optimal transport problem with an inequality constraint; their
algorithm penalizes the optimization problem in its dual representation.

In this paper, we study optimal transport by semi-martingales as introduced in \cite{tan2013optimal} and we use deep learning to estimate the optimal
drift and diffusion coefficients. 

In particular, we propose two neural network-based algorithms in the paper.
\begin{itemize}
\item In the first one, we relax the terminal constraint by
adding a penalty term to the loss function and solve it with a deep neural network. \footnote{This relaxation can be useful in the case where, depending on the constraints we put on the stochastic evolution, not all distributions are attainable, think for example of the case where we impose the process to be a martingale}
\item In the second one, we introduce the dual formulation of the problem and express it as a saddle point problem. In this case,
we utilize adversarial networks to solve it.
\end{itemize}
These methods can be widely applied to solving optimal transport as well as stochastic optimal control problems. The two methods do not need spatial discretization and hence can be potentially used for high-dimensional problems.
We illustrate our method with an application in finance, where we implement the first algorithm to the problem of optimal portfolio selection with a prescribed terminal density studied in \citet{guo2020portfolio}.

Our paper is organized as follows. In \prettyref{sec:Problem-formulation}, we formulate the
optimal transport problem and introduce the primal
problem, adding a penalization to the expected cost function.
To solve this problem, we propose a deep neural network-based algorithm and present the corresponding numerical results in \prettyref{sec:DNN-Algorithm-primal}. 
In \prettyref{sec:GANs-Algorithm-dual}, we provide
the dual representation of the primal problem and express it as a
saddle point problem. Then we devise an algorithm using adversarial
networks. We illustrate the numerical results in subsections \ref{sec:Numerical-Results} and \ref{sec:higher-dimentional-example},
including a $10$-dimensional example. Finally, in \prettyref{sec:Application-in-Portfolio},
we implement the deep learning algorithm introduced in \prettyref{sec:DNN-Algorithm-primal} for a financial application, namely
solving the problem of steering a portfolio wealth towards a prescribed
terminal density.

\section{Problem formulation\label{sec:Problem-formulation}}

Let $\mathcal{D}$ be a Polish space equipped with its Borel $\sigma$-algebra.
We denote as $C(\mathcal{D};\mathbb{R}^d)$ the space of continuous functions
on $\mathcal{D}$ with values in $\mathbb{R}^d$, $C_{b}(\mathcal{D};\mathbb{R}^d)$ the space of bounded continuous functions and $C_{0}(\mathcal{D};\mathbb{R}^d)$ the space of continuous functions, vanishing at infinity. Let $\mathcal{P}(\mathcal{D})$ be the space of Borel probability measures on $\mathcal{D}$ with a finite second moment, and $L^1(d\mu)$ be the space of $\mu$-integrable functions. Let $\mathbb{R}^{+}$ denote non-negative real numbers, $\mathbb{S}^{d}$ denote the set of $d \times d$ symmetric matrices and $\mathbb{S}_{+}^{d} \subset \mathbb{S}^{d}$ denote the set of positive semidefinite matrices.
For convenience, we often use the notation $\mathcal{E}\coloneqq[0,1]\times\mathbb{R}^d$. 
We say that a function $\phi:\mathcal{E}\rightarrow\mathbb{R}$ belongs
to $C_{b}^{1,2}(\mathcal{E})$ if $\phi\in C_{b}(\mathcal{E})$ and
$(\partial_{t}\phi,\partial_{x}\phi,\partial_{xx}\phi)\in C_{b}(\mathcal{E};\mathbb{R},\mathbb{R}^d,\mathbb{S}^d)$.

Without loss of generality, we set the time horizon $T$ to be $1$. Let $\Omega\coloneqq(\omega\in C([0,1];\mathbb{R}^{d})),$ we denote by $\mathbb{F}=(\mathcal{F}_{t})_{t\in[0,1]}$ the filtration generated by the canonical process. The process $W$ is a $d$-dimensional standard Brownian motion on the filtered
probability space $\left(\Omega,\mathcal{F},\mathbb{F},\mathbb{P}\right)$.

The stochastic process $(X_{t})_{t\in[0,1]}$ valued in $\mathbb{R}^{d}$
solves the SDE 
\begin{alignat}{1}
dX_{t} & =B(t,X_{t})dt+A(t,X_{t})dW_{t},\label{eq:SDE}\\
X_{0} & =x_{0},
\end{alignat}
where $B:\mathcal{E}\rightarrow\mathbb{R}^{d}$ and $A:\mathcal{E}\rightarrow\mathbb{R}^{d\times d}$
is defined such that $AA^{\intercal}=\mathcal{A}$.

We denote by $\rho_{t}\coloneqq\mathbb{P}\circ X_{t}^{-1}\in\mathcal{P}(\mathbb{R}^{d})$
the distribution of $X_{t}$. 
In this problem, we are given the initial distribution of the state variable $\rho_{0}\in\mathcal{P}(\mathbb{R}^{d})$,
and a prescribed terminal distribution $\bar{\rho}_{1}\in\mathcal{P}(\mathbb{R}^{d})$.
We define a convex cost function $F: \mathbb{R}^{d}\times\mathbb{S}^{d}\rightarrow\mathbb{R}^{+}\cup\{+\infty\}$
where $F(B,\mathcal{A})=+\infty$ if $\mathcal{A}\notin\mathbb{S}_{+}^{d}$.

Given $\rho_{0}$ and suitable processes $(A_{t})_{t\in[0,1]},(B_{t})_{t\in[0,1]}$,
the   distribution of $X_1$ is $\rho_{1}\coloneqq\mathbb{P}\circ X_{1}^{-1}$. We  introduce the penalty function $C(\rho_{1},\bar{\rho}_{1})$, whose role will be to penalize the deviation of $\rho_{1}$ from the target $\bar{\rho}_{1}$.

Now we are interested in the following minimization problem: \begin{problem}
\label{prob:1}With a given initial distribution $\rho_{0}$ and a
prescribed terminal distribution $\bar{\rho}_{1}$, we want to solve the infimum of the functional
\begin{alignat}{1}
V(\rho_{0},\bar{\rho}_{1}) = \inf_{\rho, B,\mathcal{A} }\left\{ \int_{\mathcal{E}}F(B_{t},\mathcal{A}_{t})d\rho(t,x)+C(\rho_{1},\bar{\rho}_{1})\right\} ,\label{eq:AB_primal}
\end{alignat}
over all $(\rho, B, \mathcal{A}) \in \mathcal{P}(\mathbb{R}^d) \times \mathbb{R}^{d} \times \mathbb{S}^{d}$  satisfying the initial distribution 
\begin{alignat}{2}
\rho(0,x) & =\rho_{0}(x) & \quad\forall x\in\mathbb{R}^{d},\label{eq:initial marginal}
\end{alignat}
and the Fokker--Planck equation 
\begin{alignat}{2}
\partial_{t}\rho(t,x)+\nabla_{x}\cdot(B(t,x)\rho(t,x))-\frac{1}{2}\sum_{i,j}\partial_{ij}(\mathcal{A}_{i,j}(t,x)\rho(t,x)) & =0 & \quad\forall(t,x)\in\mathcal{E}.\label{eq:Fokker-Planck}
\end{alignat}
\end{problem}

Because the  objective function \eqref{eq:AB_primal} is a
trade-off between the cost function $F$ and the penalty $C$, the
optimal $\rho, \mathcal{A},B$, if it exists, will not in general ensure that $\rho_{1}=\bar{\rho}_{1}$, unless the penalty function is
\begin{equation}
C(\rho_{1},\bar{\rho}_{1})=\begin{cases}
0 & \text{ if }\rho_{1}=\bar{\rho}_{1},\\
+\infty & \;\text{if}\:\rho_{1}\neq\bar{\rho}_{1},
\end{cases}\label{eq:characteristic function}
\end{equation}
and one
recovers the ``usual'' semi-martingale optimal transport problem.

We make the following assumptions which will hold throughout the paper.
\begin{assumption} 
For $ t\in(0,1]$ the probability measure $\rho_{t}$
is absolutely continuous with respect to the Lebesgue measure. 
\end{assumption}
\begin{assumption} \label{assu:C-is-continuous.} 
The penalty function
$C(\cdot,\bar{\rho}_{1}):\mathcal{P}(\mathbb{R}^d)\rightarrow\mathbb{R}^{+}$
is convex, lower semi-continuous with respect to the weak-$\ast$ convergence of $\rho_1$, and  $C(\rho_{1},\bar{\rho}_{1})=0$
if and only if $\rho_{1}=\bar{\rho}_{1}$.
\end{assumption} 

\begin{assumption} $\,$ \label{assu:finiteness of F} 
\begin{enumerate}[label=(\roman*)]
\item For $ (t,x) \in\mathcal{E}$ the function $(B, \mathcal{A})\to F(B,\mathcal{A})$ is non-negative, lower semi-continuous and strictly convex. 
\item There exist constants $m>1$ and $k>0$ such that for all $(t,x) \in \mathcal{E}$ there holds
\[
\left|B\right|^{m}+\left|\mathcal{A}\right|^{m}\leq k\left(1+F(B,\mathcal{A})\right).
\]
\item For all $(t,x)\in\mathcal{E}$, we have 
\[
\mathbb{E}\Bigl[\int_{0}^{1}\left|B\right|+\left|\mathcal{A}\right|dt\Bigr]<\infty,
\]
where $\left|\cdot\right|$ is the $L^{1}$-norm.
\end{enumerate}
\end{assumption}

With Assumptions \ref{assu:C-is-continuous.} and \ref{assu:finiteness of F}, 
we can get the following existence and uniqueness result for the minimizer of Problem \ref{prob:1} using similar convex minimization techniques as in 
\citet[Proposition 2]{loeper2006reconstruction}.  We refer the interested readers to \citet{loeper2006reconstruction} for the detailed proof.

\begin{thm}
\label{thm:uniqueness_minimizer}
Let us define 
\[
I(\rho, B, \mathcal{A}) =  \int_{\mathcal{E}}F(B_{t},\mathcal{A}_{t})d\rho(t,x)+C(\rho_{1},\bar{\rho}_{1}).
\]
Under Assumptions \ref{assu:C-is-continuous.} and \ref{assu:finiteness of F}, if $I(\rho, B, \mathcal{A}) $ is finite, then there exists a unique minimizer $(\rho, B, \mathcal{A}) \in \mathcal{P}(\mathbb{R}^d) \times \mathbb{R}^{d} \times \mathbb{S}^{d}$ for Problem \ref{prob:1}  satisfying constraints \eqref{eq:initial marginal} and \eqref{eq:Fokker-Planck}.
\end{thm}

\section{Deep Neural Network based Algorithm \label{sec:DNN-Algorithm-primal}}

In this section, we devise an algorithm to solve problem \prettyref{prob:1} by using deep learning. 
Note that given processes $(B_t)_{t \in [0,1]}$, $(A_t)_{t \in [0,1]}$ and $\rho_0$,
the density process $(\rho_t)_{t \in (0,1]}$ is fully determined (up to suitable regulation assumptions on $B$ and $A$). 
Hence, our goal is to use neural networks to search for the optimal  $(B_t)_{t \in [0,1]}$ and $(A_t)_{t \in [0,1]}$, which minimize the objective function. 

We discretize the period $[0,1]$ into $N$ constant time steps and construct a neural network $\theta_{n}$ for each time step $n \in [0, N-1]$.
We use multilayer feedforward neural networks in our application. At each time step, we wish to approximate the drift and diffusion coefficient with the neural network, i.e., $(B_{n},A_{n}) \approx\theta_{n}(X_{n})$, where $B_{n}\in\mathbb{R}^{d}, A_{n}\in\mathbb{R}^{d\times d}$. Feedforward neural networks approximate complicated nonlinear functions by a composition of simpler functions, namely
\begin{equation*}
(B_{n},A_{n})\approx\theta_{n}(X_{n}) = g_n^J  \circ g_{n}^{J-1} \circ \cdots \circ g_{n}^{1}(X_n).
\end{equation*}
For each layer, $g_{n}^{j}$ is
\begin{equation}
g_{n}^{j}(X_n) =  \sigma^{j}\left(X_{n}\textbf{W}_{n}^{j}+\textbf{b}_{n}^{j}\right),
\end{equation}
where $\textbf{W}_{n}, \textbf{b}_{n}$ are the weight matrices and bias vectors, respectively. Here, $\sigma^{j}(\cdot)$ is a component-wise nonlinear activation function, like sigmoid, ReLU, tanh, etc. In our paper, we use Leaky ReLU for all the hidden layers. For the output layer, $\sigma^{j}(\cdot)$ is the identity function. 

We denote by  $M$ the number of Monte Carlo paths. For a particular path $m \in [1, M]$, with $B_{n}^{m},A_{n}^{m}$ and $X_n^m$, we can compute the state variable in the next time step from the dynamics 
\[
X_{n+1}^{m}=X_{n}^{m}+B_{n}^{m}\Delta t+A_{n}^{m}\Delta W_{n}^{m},\quad \forall m\in[1,M],n\in[0,N-1].
\]

At the final time step $N$, we can get $M$ samples of terminal wealth $X_{N}$ from the Monte Carlo paths. 
With these samples,  we can estimate the empirical terminal distribution $\rho_{1}$ using kernel density
estimation (KDE). In particular, we use the Gaussian kernel $K_{\mathbf{H}}(x)=(2\pi)^{-d/2}|\mathbf{H}|^{-\frac{1}{2}}e^{-\frac{1}{2}x^{\intercal}\mathbf{H}^{-1}x}$
with an appropriate bandwidth matrix $\mathbf{H} \in \mathbb{S}^d$. 
Then the terminal density is estimated as $\tilde{\rho}_{1}(x)\coloneqq\frac{1}{M}\sum_{m=1}^{M}K_{\mathbf{H}}(x-X_{N}^{m})$.
Because the kernel density estimation is not an unbiased estimator of
the true density, an error is generated from estimating $\rho_{1}$
with $\tilde{\rho}_{1}$. To address this issue, we also estimate
$\bar{\rho}_{1}$ with the same KDE and use the estimated $\tilde{\bar{\rho}}_{1}$
as the target density in the training. To be precise, we first generate
a sample $\tilde{\mathbf{x}}=(\tilde{x}_{1},\tilde{x}_{2},...,\tilde{x}_{M'})$
of size $M'$ from the target density $\bar{\rho}_{1}(x)$. Then we
estimate $\bar{\rho}_{1}(x)$ with $\tilde{\bar{\rho}}_{1}(x)\coloneqq\frac{1}{M'}\sum_{i=1}^{M'}K_{\mathbf{H}}(x-\tilde{x}_{i})$.

Finally, the objective function in \eqref{eq:AB_primal} can be naturally used as the loss function in the training, 
\[
L(\theta_{n,n\in[0,N-1]})=\mathbb{E}\Bigl[\sum_{n=0}^{N-1}F(\theta_{n}(X_{n}))\Delta t\Bigr]+C(\tilde{\rho}_{1},\tilde{\bar{\rho}}_{1}).
\]
Then the training will search for the optimal neurons $\hat{\theta}_{n,n\in[0, N-1]}$ where 
\[
\hat{\theta}_{n,n\in[0,N-1]}=\arg\inf_{\theta}\left\{ \frac{1}{M}\sum_{m=1}^{M}\biggl\{\sum_{n=0}^{N-1}F(\theta_{n}(X_{n}^{m}))\Delta t\biggr\}+C(\tilde{\rho}_{1},\tilde{\bar{\rho}}_{1})\right\} .
\]

\begin{figure}[H]
\centering 
\includegraphics[scale=0.4]{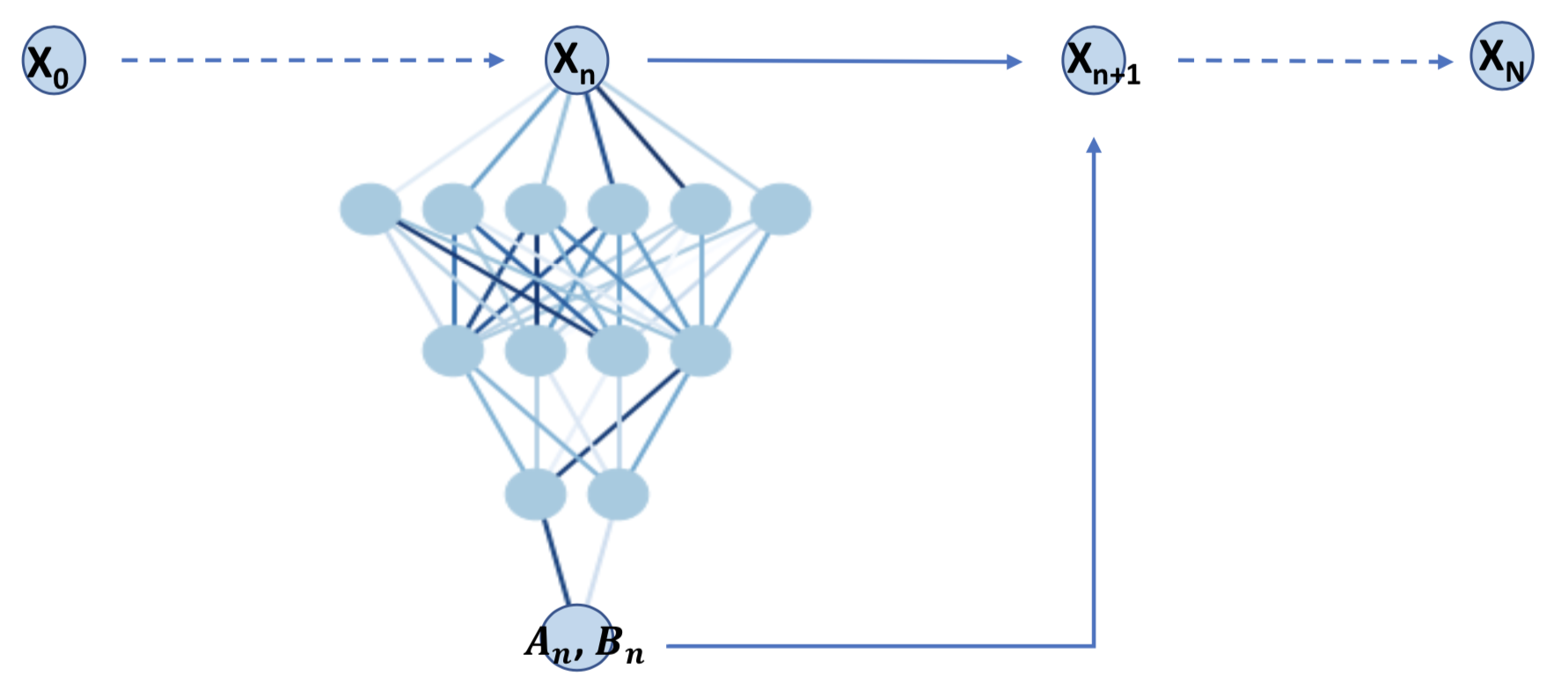} 
\caption{Structure of the DNN for one path}
\label{fig:primal_flowchart_AB} 
\end{figure}

We depict the above process in a flowchart (Figure \ref{fig:primal_flowchart_AB}), and the complete Deep Neural Network algorithm for Problem \prettyref{prob:1} is stated in Algorithm \prettyref{alg:Deep_Learning}.

\begin{algorithm}[H]
Starting from the initial condition $X_{0}^{m}=x_{0}\,\forall m\in[1,M]$.

\For{epoch $=1:100$}{

\For{$n=0:N-1$}{

With neurons $\theta_{n}$, input $X_{n}^{m}$, output $B_{n}^{m},A_{n}^{m}=\theta_{n}(X_{n}^{m})$;

$X_{n+1}^{m}=X_{n}^{m}+B_{n}^{m}\Delta t+A_{n}^{m}\Delta W_{n}^{m}$;
} Estimate the terminal distribution of $X_{N}$ with KDE as $\tilde{\rho}_{1}$;

Estimate $\tilde{\bar{\rho}}_{1}$ from $\bar\rho_1$ with the same KDE;

Define Loss function $L=\frac{1}{M}\sum_{m=1}^{M}\biggl\{\sum_{n=0}^{N-1}F(B_{n}^{m},A_{n}^{m}(A_{n}^{m})^\intercal)\Delta t\biggr\}+C(\tilde{\rho}_{1},\tilde{\bar{\rho}}_{1})$;

Train the neurons and update $\theta_{n,n\in[0,N-1]}$;
} Get the optimal $\hat{\theta}_{n,n\in[0,N-1]}$.

\caption{Deep Neural Network for Problem \prettyref{prob:1} \label{alg:Deep_Learning}}
\end{algorithm}

\subsection{Numerical Results  
\label{sec: 2_d_primal}}

In this section,  we validate Algorithm \ref{alg:Deep_Learning} with an example where
\begin{align}
F(B, \mathcal{A}) & =\left\Vert B \right\Vert ^{2}, \\
C(\rho_{1},\bar{\rho}_{1}) & =\frac{\lambda}{2} \int_{\mathbb{R}} (\rho_{1}-\bar{\rho}_{1})^{2}dx. \label{eq:penalty_primal}
\end{align}
The parameter $\lambda$ in \eqref{eq:penalty_primal} is the intensity of the penalization. To push $\rho_1$ as close as possible to the target $\bar{\rho}_1$, we want $\lambda$ be to large and we will let $\lambda = 5000$ in this case. 
We start from an initial state 
$ x_0 = \begin{bmatrix}5.0\\
5.0
\end{bmatrix}$,
and the target distribution set to be a bivariate normal  $\bar{\rho}_{1}=\mathcal{N}\left(\begin{bmatrix}5.5\\
6.0
\end{bmatrix},\begin{bmatrix}0.25 & 0.10\\
0.10 & 0.25
\end{bmatrix}\right)$.

In the training, we construct a $5$-layer  network with neurons $[100, 80, 60, 60, 40]$. The output of the network is of size $2 \times 3$,  out of which the $2 \times 1$ column represents the drift $B_n$ and the $2 \times 2$ matrix represents the diffusion coefficient $A_n$. We use a batch size of $2000$ and the Adam optimizer with a learning rate $1 \times 10^{-4}$.  We trained the network for $100$ epochs in total. 

From the training, we get an empirical terminal density $\rho_{1}$ with mean $\begin{bmatrix}5.509\\
5.994
\end{bmatrix}$ and covariance matrix $\begin{bmatrix}0.253 & 0.098\\
0.098 & 0.246
\end{bmatrix}$. We can see that the mean and covariance of the empirical terminal density are very close to the target ones. 

\begin{figure}[H]
\centering \subfloat[ Density of $X_{1}$ from Algorithm \ref{alg:Deep_Learning}\label{fig:dens_2d_rho_T}]
{{\includegraphics[scale=0.4]{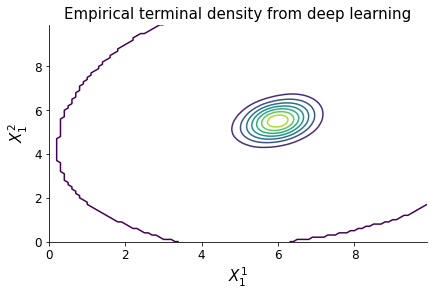} }}
\qquad{}\quad{}
\subfloat[Target Density of $X_{1}$ \label{fig:dens_2d_rho_bar}]
{{\includegraphics[scale=0.4]{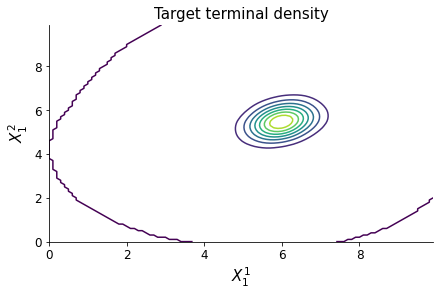} }}
\caption{Contours with a bivariate normal target distribution}
\label{fig:prime_2d}
\end{figure}

The contours of the empirical distribution and the target distribution
are presented in Figure \ref{fig:prime_2d}. 
In the later section, we will introduce a metric to measure the performance of the trained dataset.

\subsection{Merged Network}

\begin{figure}[H]
\centering \includegraphics[scale=0.25]{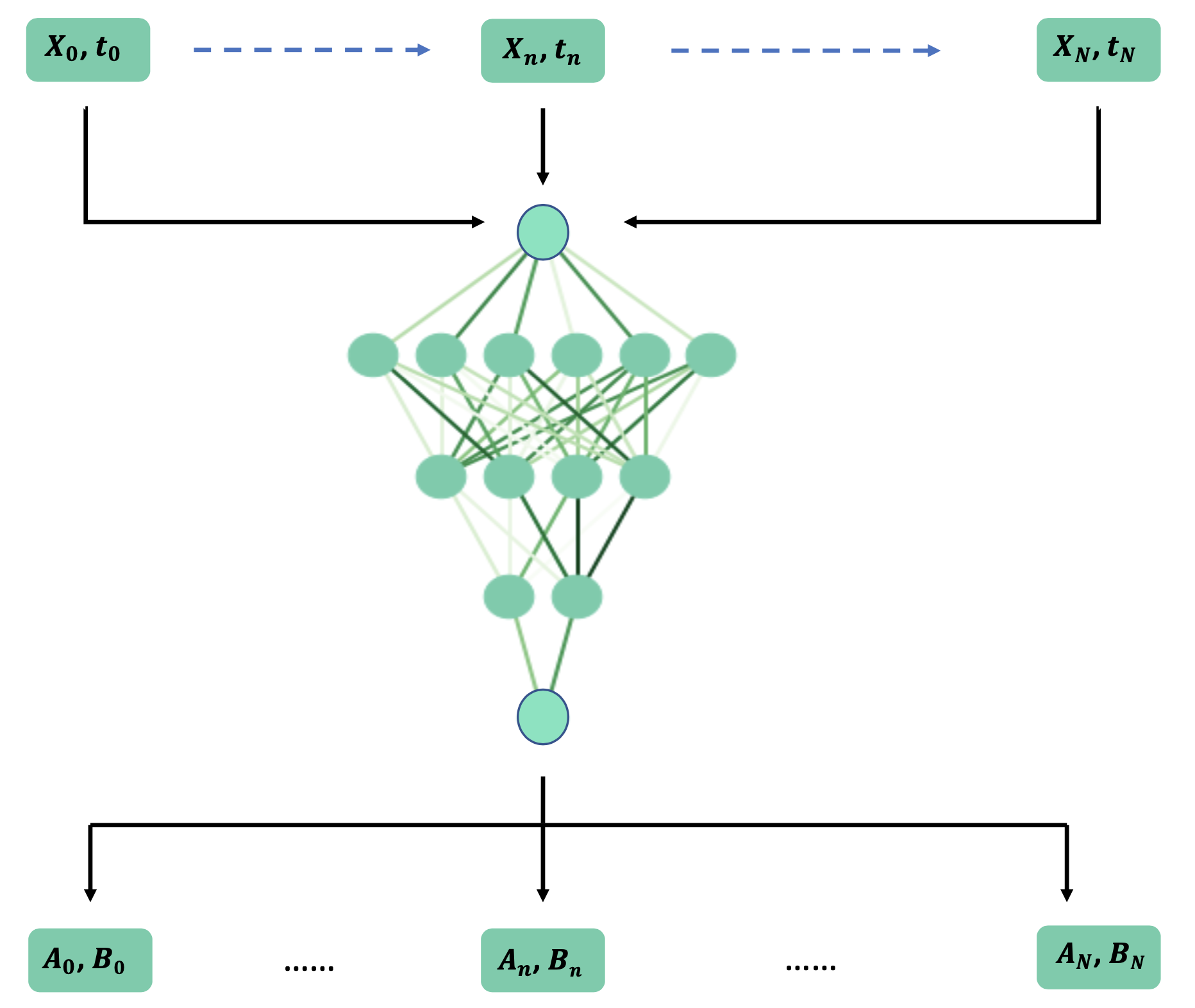} \caption{Structure of the merged network}
\label{fig:primal_single_flowchart} 
\end{figure}

The architecture described in Algorithm \ref{alg:Deep_Learning} consists
of training $N$ different feedforward neural networks (one per time
step). This architecture generates a possibly high number of weights
and bias to be estimated. Another possibility is to use one single
feedforward neural network for all time steps. At each time step, besides the state variable
$X_{n}$, we also feed the time $n$ into the network. 
We refer to this alternative architecture as a \textit{feedforward merged network}. 
The training process for this merged
network is illustrated in Figure \ref{fig:primal_single_flowchart},
and the algorithm is summarized in Algorithm \prettyref{alg:Deep_Learning_merged}.
The advantage of this architecture is that, instead of training a list of networks
$\theta_{n,n\in[0,N-1]}$, we only need to train one network $\theta$
in this algorithm, which significantly reduces the number of neurons and the complexity of the training.

The numerical experiment gives us a similar result to the one from Algorithm \ref{alg:Deep_Learning}; hence we are not presenting them repeatedly.
As we observed from the experiment, each epoch's speed of training is faster, but it will need relatively more epochs
to converge to the final value. 

\begin{algorithm}
Starting from the initial condition $X_{0}^{m}=x_{0}\,\forall m\in[1,M]$.

\For{epoch $=1:100$}{

\For{$n=0:N-1$}{

With neurons $\theta$, input $X_{n}^{m}$ and $n$, output $B_{n}^{m},A_{n}^{m}=\theta(X_{n}^{m}, n)$;

$X_{n+1}^{m}=X_{n}^{m}+B_{n}^{m}\Delta t+A_{n}^{m}\Delta W_{n}^{m}$;
} Estimate the terminal distribution of $X_{N}$ with KDE as $\tilde{\rho}_{1}$;

Estimate $\tilde{\bar{\rho}}_{1}$ from $\bar\rho_1$ with the same KDE;

Loss function $L=\frac{1}{M}\sum_{m=1}^{M}\biggl\{\sum_{n=0}^{N-1}F( B_{n}^{m}, A_{n}^{m}(A_{n}^{m})^\intercal )\Delta t\biggr\}+C(\tilde{\rho}_{1},\tilde{\bar{\rho}}_{1})$;

Train the neurons with the Adam optimizer and update $\theta$; }
Get the optimal $\hat{\theta}$.

\caption{Deep Neural Network for Problem \prettyref{prob:1} with a merged network  \label{alg:Deep_Learning_merged}}
\end{algorithm}

\section{Adversarial Network Algorithm for the dual problem\label{sec:GANs-Algorithm-dual}}

In this section, we first introduce the dual formulation of Problem \ref{prob:1} and express it as a saddle point problem. Then we propose an Adversarial Network-based algorithm to solve it, where we demonstrate high dimensional numerical examples. This algorithm is inspired by 
Generative Adversarial Networks (GANs), which were first introduced in \citet{goodfellow2014generative}.
GANs have enjoyed great empirical success in image generating and processing.
The principle behind GANs is to interpret the process of generative modeling as a competing game between two (deep) neural networks: a generator and a discriminator.
The generator network attempts to generate data that looks similar to the training data, and the discriminator network tries to identify whether the input sample is faked or true. GANs recently have attracted interests in the finance field, including simulating financial time-series data (e.g., \citealt{wiese2020quant}, \citealt{zhang2019stock}) and for asset pricing models (\citealt{chen2019deep}).

To recover the optimal transport problem, we can use an indicator function \eqref{eq:characteristic function} as the penalty functional.
Now we present the duality result for Problem  \ref{prob:1}. The following  theorem is stated  in  \citet{tan2013optimal}.
\begin{thm}
\label{thm:dual_with_indicator} When $C(\rho_{1},\bar{\rho}_{1})$
is defined as \eqref{eq:characteristic function}, there holds 
\begin{alignat}{1}
V(\rho_{0},\bar{\rho}_{1}) & =\sup_{\phi_{1}}\left\{ \int_{\mathbb{R}^{d}}\phi_{1}d\bar{\rho}_{1} - \phi_{0}d\rho_{0}\right\} ,\label{eq:classical OT dual}
\end{alignat}
where the supremum is running over all $\phi_{1}\in C_{b}^{2}(\mathbb{R}^{d})$
and $\phi_{0}$ is a viscosity solution of the Hamilton--Jacobi--Bellman equation 
\begin{gather}
\begin{cases}
-\partial_{t}\phi-\sup_{B\in\mathbb{R}^{d},\mathcal{A}\in\mathbb{S}^{d}}\Bigl[B\cdot\nabla_{x}\phi+\frac{1}{2}tr(\mathcal{A}\nabla_{x}^{2}\phi)-F(B,\mathcal{A})\Bigr]=0, & \text{in }[0,1)\times\mathbb{R}^{d},\\
\phi(1,x)=\phi_{1}(x), & \text{on }[1]\times\mathbb{R}^{d}.
\end{cases}\label{eq:HJB PDE}
\end{gather}
The function $\phi(0,x)$ can be expressed as 
\begin{equation}
\phi(0,x)=\sup_{B\in\mathbb{R}^{d},\mathcal{A}\in\mathbb{S}^{d}}\mathbb{E}\left[\phi_{1}(X_{1})-\int_{0}^{1}F(B,\mathcal{A})dt|X_{0}=x_{0}\right].\label{eq:phi_0 expression}
\end{equation}
\end{thm}

Starting from an initial state $x_{0}\in\mathbb{R}^{d}$, the initial
density $\rho_{0}(x)=\delta(x-x_{0})$, hence $\int_{\mathbb{R}^{d}}\phi_{0}d\rho_{0}=\phi_{0}(x_{0})$.
Then we can substitute the expression \eqref{eq:phi_0 expression}
into the dual form \eqref{eq:classical OT dual}, and the dual formulation
can be written as a saddle point problem: 
\begin{equation}
V(\rho_{0},\bar{\rho}_{1})=\sup_{\phi_{1}}\inf_{B\in\mathbb{R}^{d},\mathcal{A}\in\mathbb{S}^{d}}\left\{ \int_{\mathbb{R}^{d}}\phi_{1}d\bar{\rho}_{1} - \mathbb{E}\left[\phi_{1}(X_{1})-\int_{0}^{1}F(B,\mathcal{A})dt|X_{0}=x_{0}\right]\right\} .
\label{eq:zero-sum-dual}
\end{equation}

This dual saddle point formulation of the optimal transport problem is reminiscent  of GANs: GANs can be interpreted as minimax games between the generator and the discriminator, whereas our problem is a minimax game between $\phi_{1}$ and $(\mathcal{A},B)$.

Inspired by this connection, we can use Adversarial Networks
to estimate the value \eqref{eq:zero-sum-dual}. GANs is also applied by \citet{guo2019robust} to solve a robust portfolio allocation problem.
Our Adversarial Network consists of
two neural networks; one generates $\phi_{1}$ (referred to as $\phi_{1}$-generator, denoted by $\Phi$),
the other generates $A$ and $B$ (referred to as $AB$-generator, denoted by $\theta$). 
The two networks are trained iteratively: In the first phase, we train
the $AB$-generator with a loss function $L_{1}=-\mathbb{E}\left[\phi_{1}(X_{1})-\int_{0}^{1}F(B,\mathcal{A})dt|X_{0}=x_{0}\right]$.
In the second phase, given the output $X_{1}$ from the $AB$-generator,
we train the $\phi_{1}$-generator with a loss function $L_{2}= - \int_{\mathbb{R}^{d}}\phi_{1}d\bar{\rho}_{1}+\mathbb{E}\left[\phi_{1}(X_{1})-\int_{0}^{1}F(B,\mathcal{A})dt|X_{0}=x_{0}\right]$.

When the dimension $d \geq 3$, the bottleneck in this algorithm is to
compute $\int_{\mathbb{R}^{d}}\phi_{1}d\bar{\rho}_{1}$ efficiently.
To address this issue, we can use Monte Carlo integration. In particular, we
sample $\mathcal{I}$ points $\bar{x}_{i}\in\mathbb{R}^{d}$ $(i \in [1, \mathcal{I}])$
from the target distribution $\bar{\rho}_{1}$. Then we can estimate
$\int_{\mathbb{R}^{d}}\phi_{1}d\bar{\rho}_{1}$ with $\frac{1}{\mathcal{I}}\sum_{i=1}^{\mathcal{I} }\phi_{1}(\bar{x}_{i})$.
In this case, this Adversarial Network-based scheme is mesh-free.  It can now avoid the curse of dimensionality and has the potential to be applied to high dimensional problems.

\begin{figure}[H]
\centering 
\includegraphics[scale=0.3]{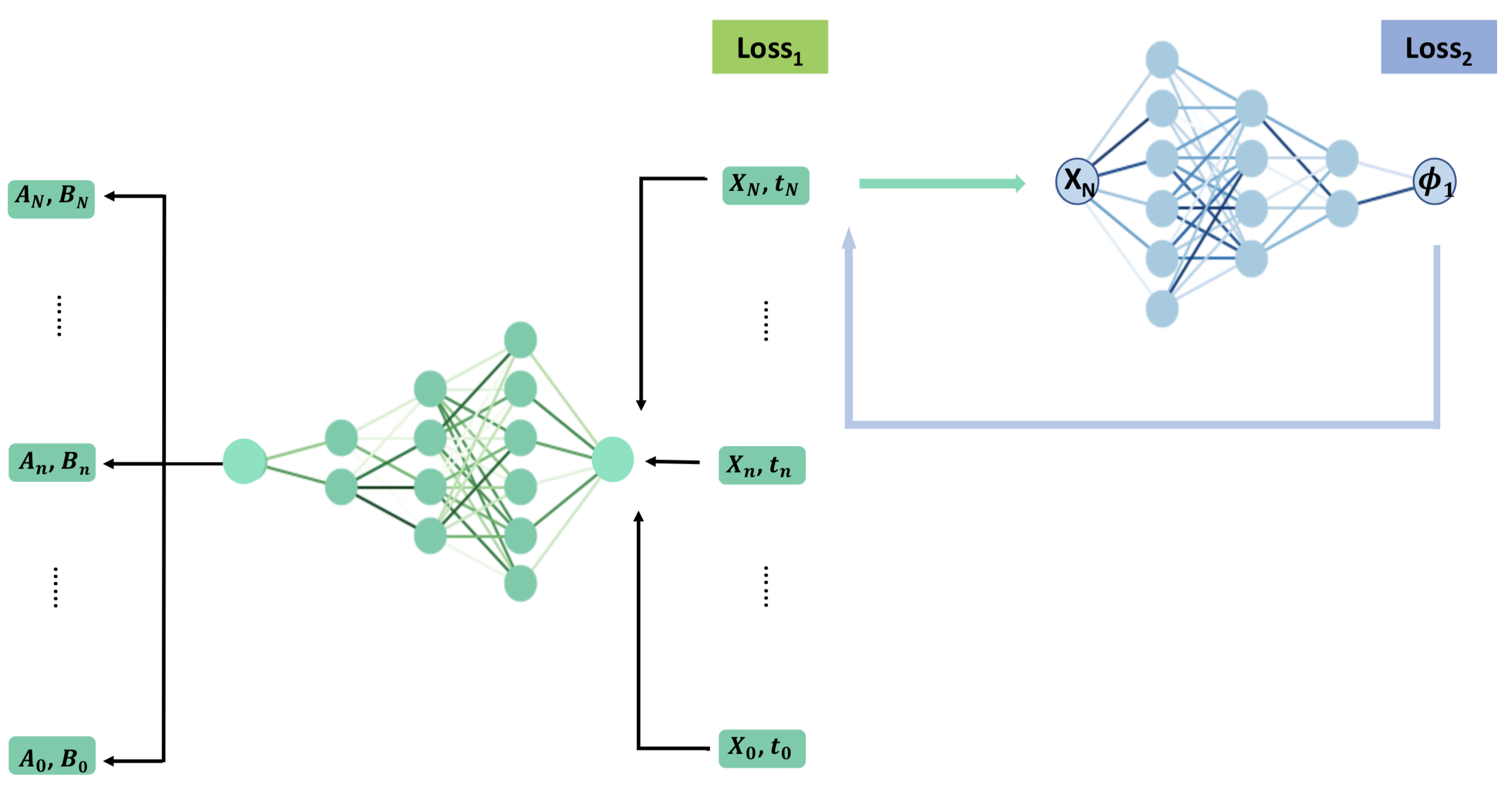} 
\caption{Adversarial Network algorithm with a merged network for the dual problem }
\label{fig:dual_gans} 
\end{figure}

We use a merged network for the $AB$-generator in this algorithm. A demonstration of the training process is illustrated in Figure \ref{fig:dual_gans}. The detailed algorithm is summarized in Algorithm \ref{alg:GANs_dual}.

\begin{algorithm}
Sample $\mathcal{I}$ points $\bar{x}_{i}\in\mathbb{R}^{d}$ $(i \in [1, \mathcal{I}])$
from the target distribution $\bar{\rho}_{1}$.

Starting from the initial condition $X_{0}^{m}=x_{0}\,\forall m\in[1,M]$:

\For{epoch $=1:\text{number of epochs}$}{ \medskip{}
 \textbf{Phase 1}: train the $AB$-generator

\For{time step $n=1:N-1$}{

With the network $\theta$, input $\{ X_{n}^m \}_{m=1}^{M}$ and
time step $n$, outputs  $(B_{n}^{m},A_{n}^{m})=\theta(X_{n}^{m},n)$;

$X_{n+1}^{m}=X_{n}^{m}+B_{n}^{m}\Delta t+A_{n}^{m}\Delta W_{n}$;
}

Loss function $L_{1}=-\frac{1}{M}\sum_{m=1}^{M}\left[\Phi(X_{N}^{m})-\sum_{n=0}^{N-1}F(B_{n}^{m},A_{n}^{m}(A_{n}^{m})^{\intercal})\Delta t\right]$;

Train the neurons with the Adam optimizer and update $\theta$. \medskip{}

\textbf{Phase 2}: train the $\phi_{1}$-generator

With the network $\Phi$, input $\{X_{N}^m \}_{m=1}^{M}$, output $\phi_{1}(X_{N}^{m})=\Phi(X_{N}^{m})$;

With the same network $\Phi$, inputs 
$\bar{x}_i, \forall i \in [1, \mathcal{I}]$, outputs $\Phi(\bar{x}_i)$;

Loss function $L_{2}=- \frac{1}{\mathcal{I}}\sum_{i=1}^{\mathcal{I} }\Phi (\bar{x}_{i}) + \frac{1}{M}\sum_{m=1}^{M}\left[\Phi(X_{N}^{m})-\sum_{n=0}^{N-1}F(B_{n}^{m},A_{n}^{m}(A_{n}^{m})^{\intercal})\Delta t\right]$;

Train the neurons with the Adam optimizer and update $\Phi$.

}

\caption{\label{alg:GANs_dual} Adversarial Network algorithm with a merged network for the
dual problem }
\end{algorithm}

\subsection{Numerical Results\label{sec:Numerical-Results}}

First, we start with a one-dimensional toy example to assess the quality of Algorithm \ref{alg:GANs_dual}. We choose a target
distribution $\bar{\rho}_{1}=\mathcal{N}(6,1)$, a cost function $F(B,\mathcal{A})=(\mathcal{A} - 0.1)^2$ and solve for $V(\rho_{0},\bar{\rho}_{1})=\sup_{\phi_{1}}\left\{ \int_{\mathbb{R}}\phi_{1}d\bar{\rho}_{1}-\phi_{0}d\rho_{0}\right\} $.

\begin{figure}[H]
\centering \subfloat[ Distribution of $X_{1}$ learnt by the Adversarial Network \label{fig:6_1_gans_inf_dens}]
{{\includegraphics[scale=0.4]{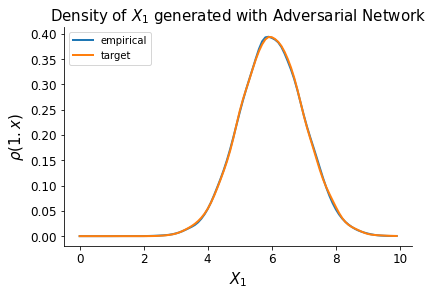} }}
\qquad{}\quad{}
\subfloat[Loss function \label{fig:6_1_gans_inf_loss}]
{{\includegraphics[scale=0.4]{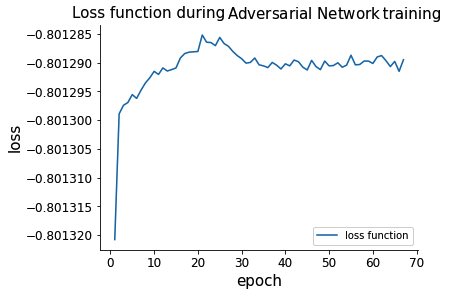} }}
\caption{Dimension $d = 1$, $\bar{\rho}_{1}=\mathcal{N}(6,1)$}
\label{fig:6_1_gans_inf}
\end{figure}

In Figure \ref{fig:6_1_gans_inf_dens}, we plot the target distribution and the distribution of $X_1$ learnt by the Adversarial Networks; in Figure \ref{fig:6_1_gans_inf_loss},  we show the corresponding loss function during the training.  
We can see that our algorithm works well in terms of attaining the target density.

Next, we demonstrate a two-dimensional example. We let $X_{n}\in\mathbb{R}^{2}$, and
we generate the drift $B_{n}\in\mathbb{R}^{2}$ and the diffusion coefficient $A_{n}\in\mathbb{R}^{2\times2}$ with the deep neural network.
We again start from an initial state 
$ x_0 = \begin{bmatrix}5.0\\
5.0
\end{bmatrix}$,
and let the cost function $F(B, \mathcal{A})=\left\Vert B\right\Vert ^{2} $ and  the target distribution set to be a bivariate normal  $\bar{\rho}_{1}=\mathcal{N}\left(\begin{bmatrix}5.5\\
6.0
\end{bmatrix},\begin{bmatrix}0.25 & 0.10\\
0.10 & 0.25
\end{bmatrix}\right)$.

We construct two networks for $(A, B)$ and $\phi_1$, respectively. The $AB$-generator has $4$ layers with neurons $[40, 30, 20, 10]$ and the $\phi_{1}$-generator has $4$ layers with neurons $[80, 60, 40, 40]$. We use a mini-batch SGD algorithm with a batch size of $1000$ and Adam optimizer with learning rate $1 \times 10^{-4}$.

\begin{figure}[H]
\centering \subfloat[Contour of $\rho_1$ from the Adversarial Network \label{fig:gans_2d_onenet_rho_T}]
{{\includegraphics[scale=0.5]{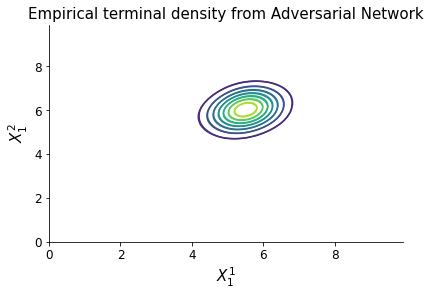} }}
\subfloat[Contour of $\bar{\rho}_1$ \label{fig:gans_2d_onenet_rho_bar}]
{{\includegraphics[scale=0.5]{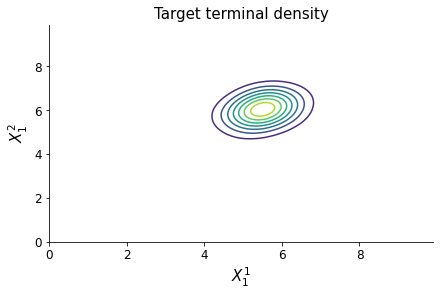} }}
\caption{Dimension $d = 2$, bivariate normal target distribution}
\label{fig:dual_2d}
\end{figure}

The empirical terminal density $\rho_{1}$ after the training has
mean $\begin{bmatrix}5.497\\
6.007
\end{bmatrix}$ and covariance matrix $\begin{bmatrix}0.253 & 0.099\\
0.099 & 0.251
\end{bmatrix}$, which are very close to the target mean and covariance. 
The contours of the empirical and target distributions are shown in Figure \ref{fig:dual_2d}. This result is similar to the one we got from Algorithm \prettyref{alg:Deep_Learning}.

\subsection{Higher dimensional examples \label{sec:higher-dimentional-example}}

In this section, we apply the Adversarial Network-based algorithm to high dimensional examples.
When the dimension $d\geq3$, Algorithm \ref{alg:Deep_Learning} will be less effective. 
It is because the kernel density estimation requires a spatial discretization, and the number of grid will increase exponentially with the dimension. The estimation becomes less accurate, and the computing cost  will be expensive in this way. 
However, the Adversarial Network-based algorithm for the dual problem is free from spatial discretization. Hence it is
less affected by an increase in dimension. 

In the following, we use a cost function $F(B, \mathcal{A})=\left\Vert B\right\Vert ^{2}+\left\Vert \mathcal{A}\right\Vert ^{2}$
and solve the dual problem \eqref{eq:zero-sum-dual} with $d=5$ and $d=10$, respectively.

In the $5-$d example, we set 
$ x_0 = \begin{bmatrix}5.0\\
5.0\\
5.0\\
5.0\\
5.0
\end{bmatrix}$, and let $\bar{\rho}_{1}$ be a multivariate normal
distribution, where $\bar{\rho}_{1}=\mathcal{N}\left(\begin{bmatrix}5.5\\
6.0\\
5.8\\
6.0\\
6.2
\end{bmatrix},\begin{bmatrix}0.25 & 0.10 & 0.10 & 0.10 & 0.10\\
0.10 & 0.25 & 0.10 & 0.10 & 0.10\\
0.10 & 0.10 & 0.25 & 0.10 & 0.10\\
0.10 & 0.10 & 0.10 & 0.25 & 0.10\\
0.10 & 0.10 & 0.10 & 0.10 & 0.25
\end{bmatrix}\right)$ \footnote{For short, we represent this prescribed distribution as $\bar{\rho}_1 = \mathcal{N}(\bar{\mu}, \bar{\Sigma})$.}.

Because the dimension is higher and the cost function is more complicated in this case, we need a bigger and deeper network. The $AB$-generator now has  $5$ layers with neurons $[400, 300, 200, 200, 150]$, the $\phi_1$-generator keeps the same as before. 
Using 20,000 out-of-sample points of $X_{1}$ from the trained model,
the empirical terminal density $\rho_{1}$ from Algorithm \ref{alg:GANs_dual} has mean $\begin{bmatrix}5.4929\\
5.9749\\
5.7855\\
5.9825\\
6.1944
\end{bmatrix}$ and covariance matrix $\begin{bmatrix}0.2515 & 0.1041 & 0.1038 & 0.1037 & 0.1011\\
0.1041 & 0.2372 & 0.0926 & 0.0937 & 0.0885\\
0.1038 & 0.0926 & 0.2541 & 0.1103 & 0.0979\\
0.1037 & 0.0937 & 0.1103 & 0.2686 & 0.0906\\
0.1011 & 0.0885 & 0.0979 & 0.0906 & 0.2475
\end{bmatrix}$.

\begin{figure}[H]
\centering
\subfloat[margin 1]{\label{}\includegraphics[width=.47\linewidth]{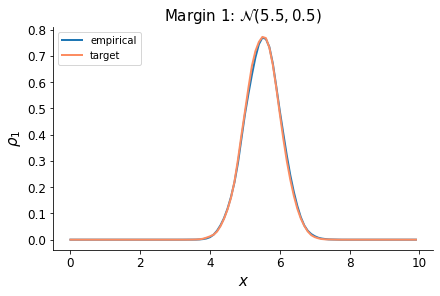}}\hfill
\subfloat[margin 2]{\label{}\includegraphics[width=.47\linewidth]{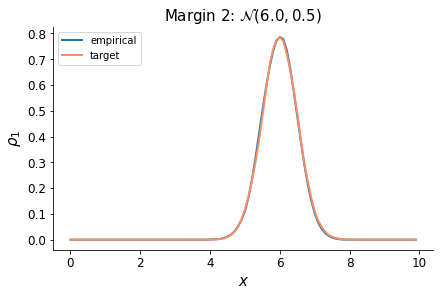}}\par
\subfloat[margin 3]{\label{}\includegraphics[width=.47\linewidth]{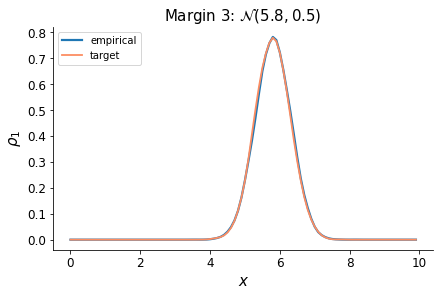}}\hfill
\subfloat[margin 4]{\label{}\includegraphics[width=.47\linewidth]{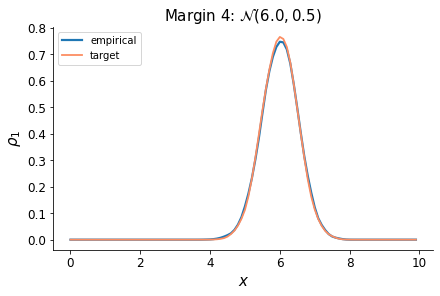}}
\par
\subfloat[margin 5]{\label{}\includegraphics[width=.47\linewidth]{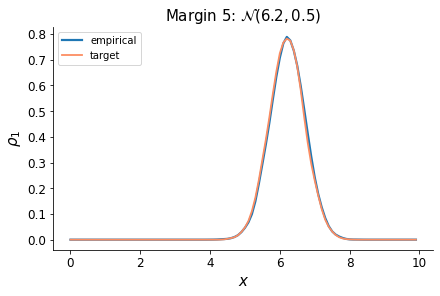}}
\caption{Marginal distributions of the 5-d sample}
\label{fig:5_d_marginal}
\end{figure}

\begin{figure}[H]
\centering
\subfloat[margin 1]{\label{}\includegraphics[width=.47\linewidth]{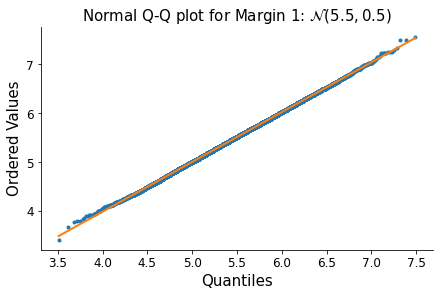}}\hfill
\subfloat[margin 2]{\label{}\includegraphics[width=.47\linewidth]{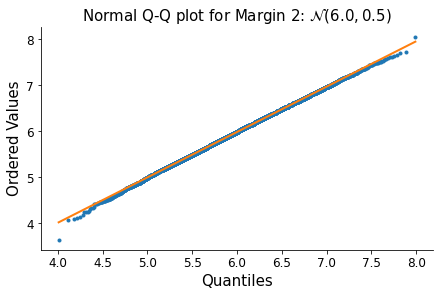}}\par
\subfloat[margin 3]{\label{}\includegraphics[width=.47\linewidth]{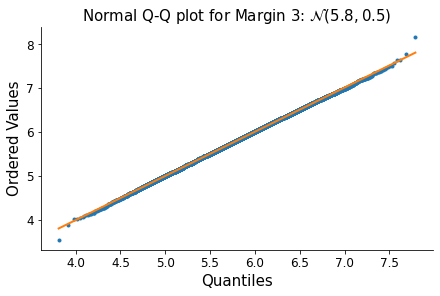}}\hfill
\subfloat[margin 4]{\label{}\includegraphics[width=.47\linewidth]{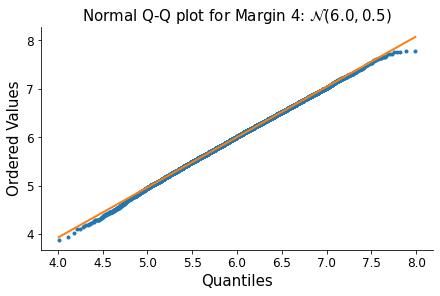}}
\par
\subfloat[margin 5]{\label{}\includegraphics[width=.47\linewidth]{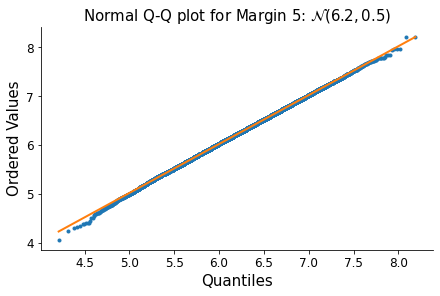}}
\caption{Normal Q-Q plot  for margins}
\label{fig:5_d_marginal_qq}
\end{figure}

We use graphical tools to help us assess if the data plausibly come from the prescribed multivariate normal distribution.
First, we compare the marginal distributions of $\rho_1$ with the ones of $\bar{\rho}_1$ in Figure \ref{fig:5_d_marginal}. We can see that the empirical marginal distributions are  consistent with the theoretical ones for all of the $5$ marginals. 
Furthermore, we make a Q-Q plot for each margin, where we plot the quantiles of the empirical marginal distribution against the quantiles of the theoretical normal distribution. As we can see in Figure \ref{fig:5_d_marginal_qq}, for every margin, the points form a roughly straight line. There are a few outliers at the two ends of the Q-Q plots, which is to be expected as Monte Carlo methods converge slowly for extreme values.
Overall, the above two figures, in addition to the mean and covariance matrix, can verify the assumption that our dataset follows the target distribution $\bar{\rho}_1$ closely.

\begin{figure}[H]
\centering
\subfloat[]{\label{}\includegraphics[width=.47\linewidth]{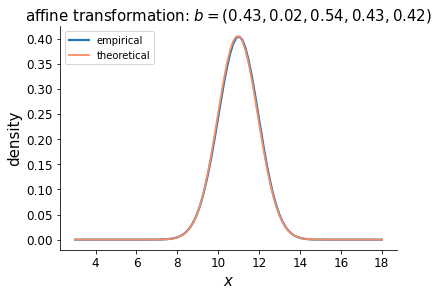}}\hfill
\subfloat[]{\label{}\includegraphics[width=.47\linewidth]{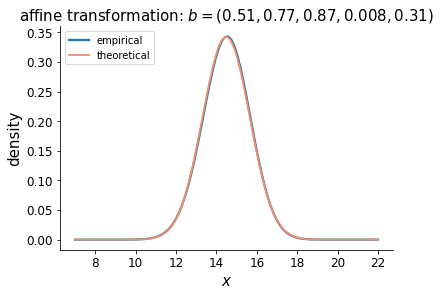}}\par
\subfloat[]{\label{}\includegraphics[width=.47\linewidth]{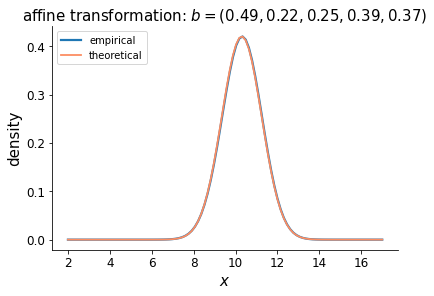}}\hfill
\subfloat[]{\label{}\includegraphics[width=.47\linewidth]{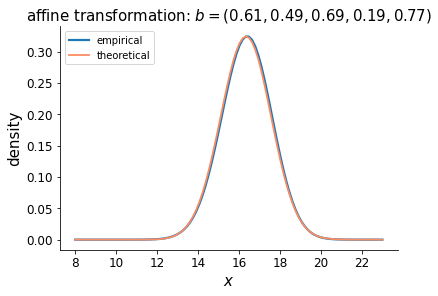}}
\par
\subfloat[]{\label{}\includegraphics[width=.47\linewidth]{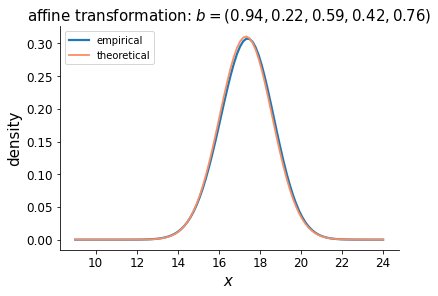}}
\caption{Distributions after affine transformations}
\label{fig: 5_d_affine}
\end{figure}

\begin{figure}[H]
\centering
\subfloat[]{\label{}\includegraphics[width=.47\linewidth]{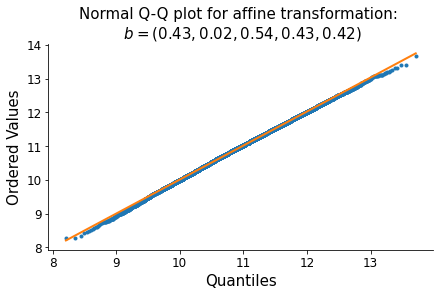}}\hfill
\subfloat[]{\label{}\includegraphics[width=.47\linewidth]{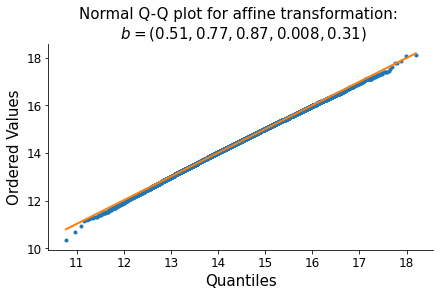}}\par
\subfloat[]{\label{}\includegraphics[width=.47\linewidth]{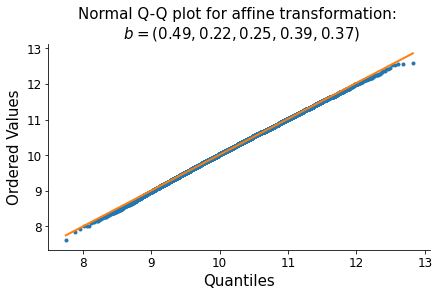}}\hfill 
\subfloat[]{\label{}\includegraphics[width=.47\linewidth]{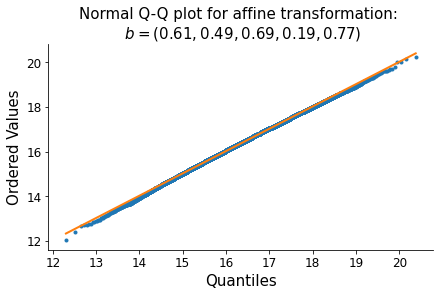}}\par
\subfloat[]{\label{}\includegraphics[width=.47\linewidth]{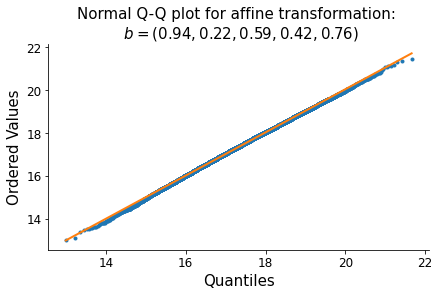}}
\caption{Normal Q-Q plots after affine transformations}
\label{fig:5_d_affine_qq}
\end{figure}

Next, we project the dataset $X_1$ on a random direction $b \in \mathbb{R}^5$. For a multivariate Gaussian distribution, the distribution after an affine transformation is univariate Gaussian.
Let $X_{\text{affine}} = X_1 b$. After the transformation, $X_{\text{affine}}$ should follow a normal distribution with mean $\mu_{\text{affine}} = \bar{\mu}^{\intercal} b$ and variance $\sigma^2_{\text{affine}} = b^{\intercal} \bar{\Sigma} b$.
In the following figures, we illustrate the results of five random affine transformations. We compare the distributions of $X_{\text{affine}} $ with the theoretical distributions (i.e., $\mathcal{N}(\mu_{\text{affine}} ,\sigma^2_{\text{affine}})$) in Figure \ref{fig: 5_d_affine}, and present their Q-Q plots in Figure \ref{fig:5_d_affine_qq}. From these two figures, we can see that the dataset after the affine transformation follows the theoretical distribution closely.

\begin{figure}[H]
\centering 
\includegraphics[scale=0.5]{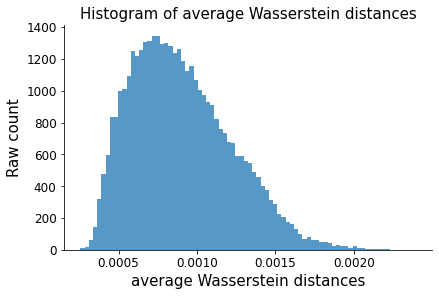} 
\caption{Average Wasserstein distances for the affine transformations of the 5 dimensional samples}
\label{fig:5d_Wassertein}
\end{figure}

Here, we introduce a loss metric to measure the quality of the dataset $X_1$.
We will apply  40,000 random affine transformations on the 5-dimensional dataset we got from the Adversarial Network.  
For each affine transformation, we use the Wasserstein distance to measure the difference between the empirical distribution and the theoretical distribution after the transformation. 
The Wasserstein distance $W_p(P,Q)$ arises from the idea of optimal transport: intuitively, it measures how far you have to move the mass of $P$ to turn it into $Q$. 
In particular, let $x_1, . . . , x_n$ be an ordered dataset from the empirical distribution $P$, and $y_1, . . . , y_n$ be an ordered dataset of the same size from the distribution $Q$, then the distance takes a very simple function of the order statistics:
$W_2(P,Q) = \left( \sum_{i=1}^{n}  ||x_i - y_i||^2  \right) ^{\frac{1}{2}}$. 
Inspired by this, we use a variation form of the 2-Wasserstein distance as the loss metric: for the $k$-th affine transformation,
define the \textit{average Wasserstein distance} between the empirical distribution and the theoretical distribution as
$W(empirical, theoretical) \coloneqq  \frac{1}{n } \sum_{i=1}^{n} \frac{1}{\bar{\sigma}_{k}^2} ||x_i - y_i||^2, \forall k \in [1, K] $, where $\bar{\sigma}_{k}$ is the standard deviation of the $k-$th theoretical distribution.
In this case, we have $n =$ 20,000 and $K =$40,000. 
To be specific, for the $k$-th transformation, we get the 20,000 ordered sample points of $X_1$ and 20,000 quantiles from the theoretical distribution, then we compute the \textit{average Wasserstein distance} between the empirical and theoretical distribution. 
In Figure \ref{fig:5d_Wassertein}, we show the histogram of the 40,000 average Wasserstein distances after affine transformations. 
 
\begin{figure}[H]
\centering 
\includegraphics[scale=0.5]{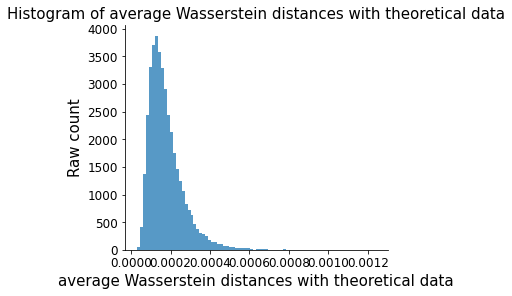} 
\caption{Average Wasserstein distances for the samples generated from theoretical distributions}
\label{fig:5d_Wassertein_theo}
\end{figure}

As a comparison, after each affine transformation, we also generate 20,000 random points from the  theoretical distribution. Then we compute the average Wasserstein distances using the 20,000 quantiles from the theoretical distribution and the ordered samples generated from the theoretical distribution.
We also present this histogram in Figure \ref{fig:5d_Wassertein_theo}. 
We can use this `correct answer' to evaluate the performance in Figure \ref{fig:5d_Wassertein}. 

Finally, we can go further to a $10$-dimensional example, where we set 

$x_ 0 = \begin{bmatrix}5.0, & 5.0,& 5.0, & 5.0, & 5.0, & 5.0, & 5.0, & 5.0, & 5.0, & 5.0 \end{bmatrix}^{\intercal}$ and the terminal distribution 
\[
\bar{\rho}_{1}=\mathcal{N}\left(\begin{bmatrix}5.5\\
6.0\\
5.8\\
6.0\\
6.2\\
5.5\\
6.0\\
5.8\\
6.0\\
6.2
\end{bmatrix},\begin{bmatrix}0.25 & 0.10 & 0.10 & 0.10 & 0.10 & 0.10 & 0.10 & 0.10 & 0.10 & 0.10\\
0.10 & 0.25 & 0.10 & 0.10 & 0.10 & 0.10 & 0.10 & 0.10 & 0.10 & 0.10\\
0.10 & 0.10 & 0.25 & 0.10 & 0.10 & 0.10 & 0.10 & 0.10 & 0.10 & 0.10\\
0.10 & 0.10 & 0.10 & 0.25 & 0.10 & 0.10 & 0.10 & 0.10 & 0.10 & 0.10\\
0.10 & 0.10 & 0.10 & 0.10 & 0.25 & 0.10 & 0.10 & 0.10 & 0.10 & 0.10\\
0.10 & 0.10 & 0.10 & 0.10 & 0.10 & 0.25 & 0.10 & 0.10 & 0.10 & 0.10\\
0.10 & 0.10 & 0.10 & 0.10 & 0.10 & 0.10 & 0.25 & 0.10 & 0.10 & 0.10\\
0.10 & 0.10 & 0.10 & 0.10 & 0.10 & 0.10 & 0.10 & 0.25 & 0.10 & 0.10\\
0.10 & 0.10 & 0.10 & 0.10 & 0.10 & 0.10 & 0.10 & 0.10 & 0.25 & 0.10\\
0.10 & 0.10 & 0.10 & 0.10 & 0.10 & 0.10 & 0.10 & 0.10 & 0.10 & 0.25
\end{bmatrix}\right).
\]

Using the same network constructed in the $5$-d example, 
the empirical terminal density $\rho_{1}$ learnt from the Adversarial Network has mean

$\begin{bmatrix}5.4711, &  5.9905, &  5.8023, &  5.9680, &  6.2131, & 5.4297, &  5.9932, &
 5.7773, & 5.9596, & 6.1853 \end{bmatrix}^{\intercal}$

and covariance matrix

$\begin{bmatrix}
0.2522 & 0.1155  & 0.1891  & 0.1305  & 0.1288   & 0.1180   & 0.1016   &  0.1148  & 0.1114 &  0.1217 \\
0.1155 &  0.2494 & 0.1346 & 0.1275 & 0.0747   & 0.1170   & 0.1265   &  0.0837 & 0.1217    &  0.1042\\
0.1891  & 0.1346  & 0.2763 & 0.1140 & 0.0939  & 0.1525   & 0.1211    &  0.1417   &  0.1298 &  0.0857\\
0.1305 & 0.1275  & 0.1140  & 0.2622  & 0.1271   & 0.1163    &  0.0892 &  0.1094  &  0.1113 &  0.1067\\
0.1288 & 0.0747  & 0.0939 & 0.1271  & 0.1858  & 0.1325   & 0.0901   &  0.0840 &  0.1120 &  0.0589\\
0.1180  & 0.1170  & 0.1525  & 0.1163  & 0.1325   & 0.2696   &  0.1435 &  0.1121  &  0.1251 &  0.1178\\
0.1016 & 0.1265   & 0.1211   &  0.0892 & 0.0901  & 0.1435   &   0.2578 &  0.1111   &  0.0751 &  0.0843\\
0.1148  & 0.0837 & 0.1417  & 0.1094  & 0.0840 & 0.1121    &  0.1111    & 0.2548  & 0.1076 &  0.1185\\
0.1114  & 0.1217   & 0.1298  & 0.1113  & 0.1120    &  0.1251   & 0.0751  & 0.1076  &   0.2455 &  0.0895\\
0.1217 & 0.1042   & 0.0857  & 0.1067 & 0.0589 & 0.1178     &  0.0843 & 0.1185  & 0.0895 &  0.2368
\end{bmatrix}$.

\begin{figure}[H]
\centering
\subfloat[margin 1]{\label{}\includegraphics[width=.38\linewidth]{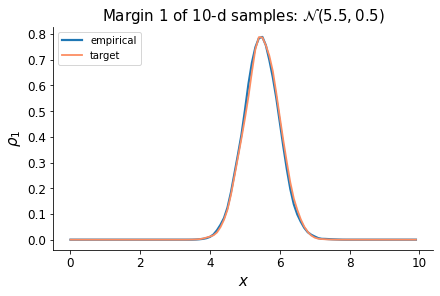}}\hfill
\subfloat[margin 2]{\label{}\includegraphics[width=.38\linewidth]{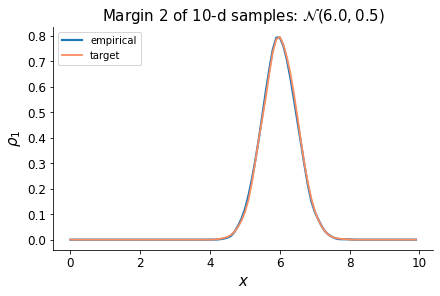}}\par
\subfloat[margin 3]{\label{}\includegraphics[width=.38\linewidth]{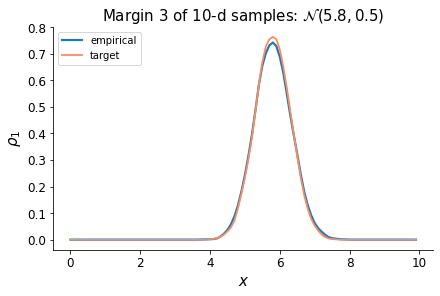}}\hfill
\subfloat[margin 4]{\label{}\includegraphics[width=.38\linewidth]{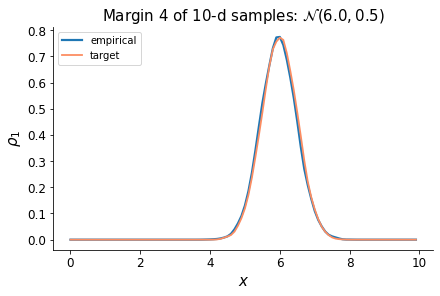}}\par
\subfloat[margin 5]{\label{}\includegraphics[width=.38\linewidth]{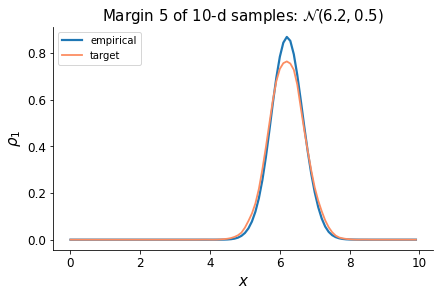}}\hfill
\subfloat[margin 6]{\label{}\includegraphics[width=.38\linewidth]{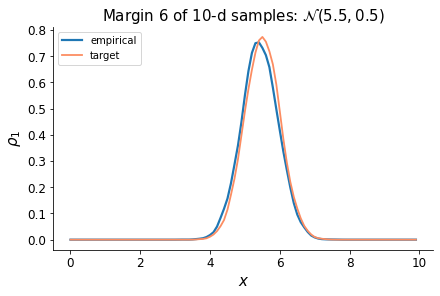}}\par
\subfloat[margin 7]{\label{}\includegraphics[width=.38\linewidth]{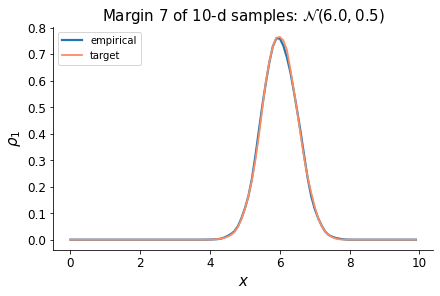}}\hfill
\subfloat[margin 8]{\label{}\includegraphics[width=.38\linewidth]{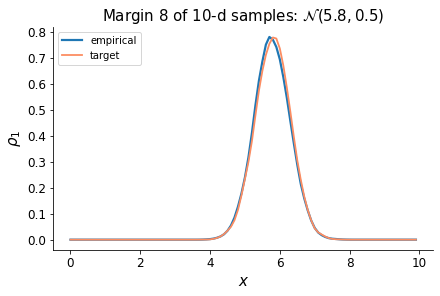}}\par
\subfloat[margin 9]{\label{}\includegraphics[width=.38\linewidth]{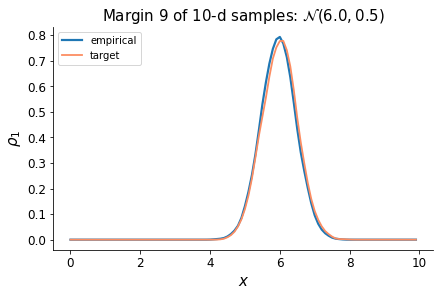}}\hfill
\subfloat[margin 10]{\label{}\includegraphics[width=.38\linewidth]{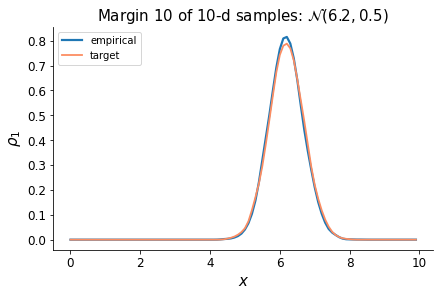}}
\caption{Marginal distributions for the 10-d sample}
\label{fig: 10_d_margins}
\end{figure}

\begin{figure}[H]
\centering
\subfloat[margin 1]{\label{}\includegraphics[width=.38\linewidth]{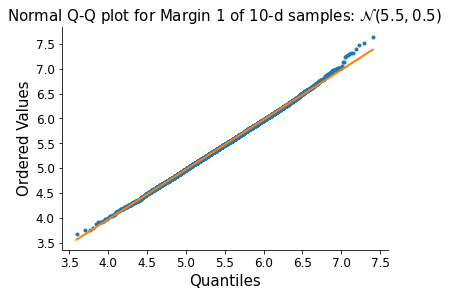}}\hfill
\subfloat[margin 2]{\label{}\includegraphics[width=.38\linewidth]{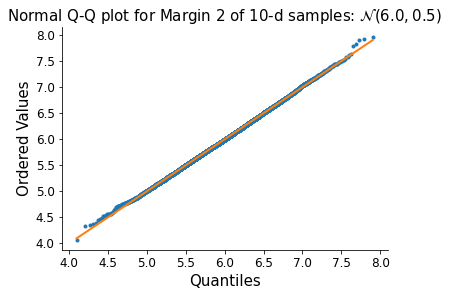}}\par
\subfloat[margin 3]{\label{}\includegraphics[width=.38\linewidth]{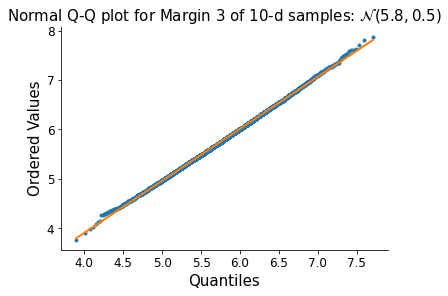}}\hfill
\subfloat[margin 4]{\label{}\includegraphics[width=.38\linewidth]{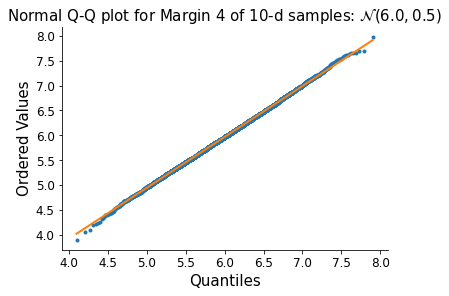}}\par
\subfloat[margin 5]{\label{}\includegraphics[width=.38\linewidth]{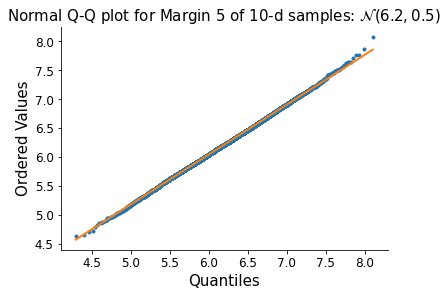}}\hfill
\subfloat[margin 6]{\label{}\includegraphics[width=.38\linewidth]{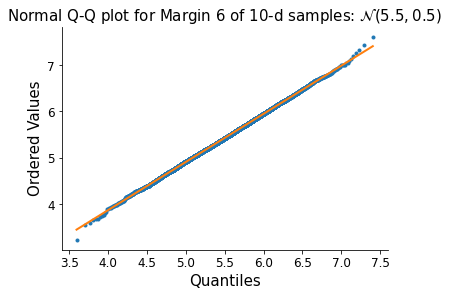}}\par
\subfloat[margin 7]{\label{}\includegraphics[width=.38\linewidth]{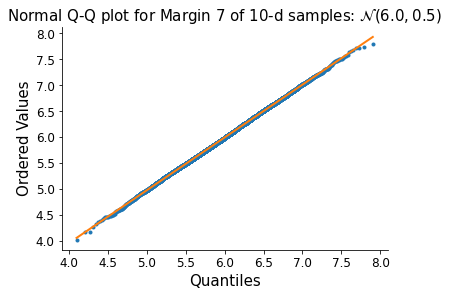}}\hfill
\subfloat[margin 8]{\label{}\includegraphics[width=.38\linewidth]{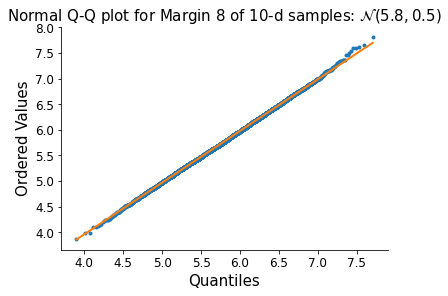}}\par
\subfloat[margin 9]{\label{}\includegraphics[width=.38\linewidth]{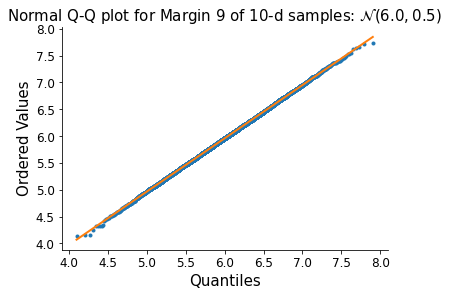}}\hfill
\subfloat[margin 10]{\label{}\includegraphics[width=.38\linewidth]{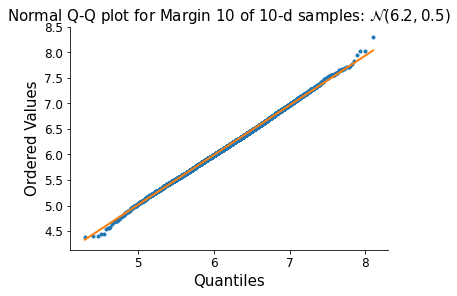}}
\caption{Normal Q-Q plots for margins of the 10-d sample}
\label{fig:10_d_margins_qq}
\end{figure}

We use a similar way to check the multivariate normality of the empirical distribution. In Figure \ref{fig: 10_d_margins}, we plot the marginal distributions of the sample and the target marginal distributions. For most of the margins, the empirical distributions are consistent with the target ones. For margin $5$, we can see the empirical one has a narrower shape (smaller variance), and the empirical margin $6$ is shifted to the left (smaller mean). These imperfectness can also be seen in the mean and covariance matrix of $\rho_1$. Then we present the Q-Q plots of the $10$ marginal distributions in Figure \ref{fig:10_d_margins_qq}. Except some outliers at the two ends, the scatters are roughly straight. The Q-Q plots verify that the margins are all normally distributed. 

\begin{figure}[H]
\centering 
\includegraphics[scale=0.5]{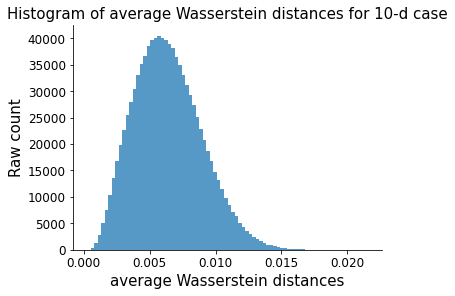} 
\caption{Average Wasserstein distances for the transformed 10-dimensional samples}
\label{fig:10d_Wassertein}
\end{figure}

\begin{figure}[H]
\centering 
\includegraphics[scale=0.5]{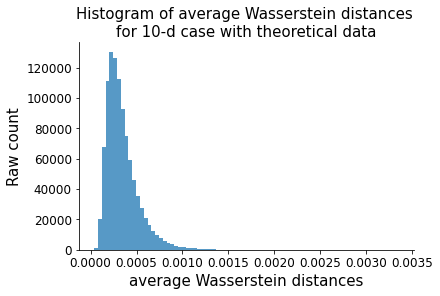} 
\caption{Average Wasserstein distances for the samples generated from theoretical distributions}
\label{fig:10d_Wassertein_theo}
\end{figure}

Then, we apply  1,000,000 affine transformations on the 10-dimensional dataset we got from the Adversarial Network. Again, we use the previously defined \textit{average Wasserstein distance} as a metric to quantify the difference between empirical distribution and the theoretical distribution after the affine transformation. In Figure \ref{fig:10d_Wassertein} we plot the histogram of the 1,000,000 average Wasserstein distances for the transformed 10-d samples, and in Figure \ref{fig:10d_Wassertein_theo} we illustrate the average Wasserstein distances for the samples generated from theoretical distributions.

We want to emphasize that, compared to other research where the prescribed distributions are restricted to Gaussian, our methods apply to a general choice of target distributions $\bar{\rho}_1$, such as heavy-tailed and asymmetric distributions. We choose multivariate normal target distributions only because the result is easier to be verified. In the following sections, we will demonstrate examples where $\bar{\rho}_1$ is not Gaussian.

\section{Application to portfolio allocation\label{sec:Application-in-Portfolio}}

Now we are applying the deep learning algorithm (Algorithm \prettyref{alg:Deep_Learning}) to the portfolio allocation problem introduced in \citet{guo2020portfolio}, where the goal is to reach a prescribed wealth distribution $\bar{\rho}_1$ at the final time from an initial wealth $x_0$. We study this problem with the tools of optimal mass transport.

In this particular application, we consider a portfolio with $d$ risky assets and one risk-free asset,
the risk-free interest $r$ being set to 0 for simplicity. We assume the drift $\mu_t:\mathcal{E} \rightarrow \mathbb{R}^{d}$ and covariance matrix $\Sigma_t: \mathcal{E} \rightarrow \mathbb{S}_{+}^{d}$ of the risky assets are known Markovian processes. 
The price process of the risky assets is denoted by $S_{t}\in\mathbb{R}^{d}$
$(0\leq t\leq1)$, and the $i$th element of $S_{t}$ follows the
semimartingale 
\begin{equation}
\frac{dS_{t}^{i}}{S_{t}^{i}}=\mu_t^{i}dt+\sum_{j=1}^{d}\sigma_t^{ij}dW_{t}^{j},\quad1\leq i\leq d,\label{eq:dynamics of S}
\end{equation}
where $\sigma_t\coloneqq\Sigma_t^{\frac{1}{2}}\in\mathbb{R}^{d\times d}$
is the diffusion coefficient matrix.

The process $\alpha = (\alpha_t)_{t \in [0,1]}$ is a Markovian control. For $t\in[0,1]$, the portfolio allocation strategy $\alpha_{t}\in\mathbb{R}^{d}$
represents the proportion of the total wealth invested into the $d$
risky assets, and $1-\sum_{i=1}^{d}\alpha_{t}^{i}$ is the proportion
invested in the risk-free asset. We define the concept of \textit{admissible
control} as follows.
\begin{defn}
\label{def: An-admissible-alpha}An admissible control process $\alpha$
for the investor on $[0,1]$ is a progressively measurable process
with respect to $\mathbb{F}$, taking values in a compact convex set
$K\subset\mathbb{R}^{d}$. The set of all admissible $\alpha$ is
compact and convex, denoted by $\mathcal{K}$. 
\end{defn}
We denote by $X_{t}\in\mathbb{R}$ the portfolio wealth at time $t$.
Starting from an initial wealth $x_{0}$, the wealth of the self-financing
portfolio evolves as follows, 
\begin{alignat}{1}
dX_{t} & =X_{t}\alpha_{t}^{\intercal}\mu_t dt+X_{t}\alpha_{t}^{\intercal}\sigma_t dW_{t},\label{eq:SDE_application}\\
X_{0} & =x_{0}.
\end{alignat}

In this case, $\rho_{t}\coloneqq\mathbb{P}\circ X_{t}^{-1}\in\mathcal{P}(\mathbb{R})$
is the distribution of the portfolio wealth $X_{t}$. With a convex cost function $F(\alpha_{t}):K \rightarrow \mathbb{R}$,
our objective function in this application is 
\begin{alignat}{1}
\inf_{\alpha,\rho}\left\{ \int_{\mathcal{E}}F(\alpha_{t})d\rho(t,x)+C(\rho_{1},\bar{\rho}_{1})\right\} ,\label{eq:primal_application}
\end{alignat}
where the feasible $(\alpha,\rho)$ in \eqref{eq:primal_application} should satisfy
the initial distribution 
\begin{alignat}{2}
\rho(0,x) & =\rho_{0}(x) & \quad\forall x\in\mathbb{R},
\end{alignat}
and the Fokker--Planck equation 
\begin{alignat}{2}
\partial_{t}\rho(t,x)+\partial_{x}(\alpha_{t}^{\intercal}\mu_t x\rho(t,x))-\frac{1}{2}\partial_{xx}(\alpha_{t}^{\intercal}\Sigma_t\alpha_{t}x^{2}\rho(t,x)) & =0 & \quad\forall(t,x)\in\mathcal{E}.
\end{alignat}

In the formulation in \eqref{eq:SDE_application}, the drift of the state variable $B(t,x) =\alpha_{t}^{\intercal}\mu_t x $ and the diffusion $\mathcal{A}(t,x) = \alpha_{t}^{\intercal}\Sigma_t \alpha_{t} x^2$. It is easy to check that $\mathcal{A}$ and $B$ are constrained in such a way that $\mathcal{A} \geq \frac{B^{2}}{\left\Vert \nu_t \right\Vert ^{2} }$, where $\nu_t \coloneqq\Sigma_t^{-\frac{1}{2}}\mu_t$. 
Note that, due to the inequality constraint on the drift and diffusion, not all target distributions are attainable in this application. Therefore, we mainly utilize Algorithm \ref{alg:Deep_Learning} to solve this problem, it can be applied to all types of terminal distributions, attainable or unattainable.

\subsection{Numerical Results}

Now we implement Algorithm \ref{alg:Deep_Learning} to solve Problem \ref{prob:1} with various penalty functionals.
In the experiment, we use $M=512000$ Monte Carlo paths, and $N=64$ time steps. We discretize the time horizon into constant steps with a
step size $\Delta t=1/N$ and the spatial domain into $I=100$
constant grids. At each time step, the neural network $\theta_n$ consists of
3 layers, with neurons $[60, 40, 20]$. We feed the
neural network with sequential mini-batches of size $1024$ and trained
it for $100$ epochs.

\subsubsection{Squared Euclidean Distance}

In this example, we assume there is one risky asset in the portfolio, i.e., $\alpha_{t}\in\mathbb{R}$. We use the  squared Euclidean distance as the penalty functional and a cost function $F(\alpha_t) = (\alpha_{t} - 0.5)^2$.
The objective function is
\[
V(\rho_{0},\bar{\rho}_{1})=\inf_{\alpha}\left\{ \mathbb{E}\Bigl[\int_{0}^{1} 
(\alpha_{t} - 0.5)^2 dt\Bigr]+\lambda\int_{\mathbb{R}}\frac{1}{2}\left(\rho_{1}(x)-\bar{\rho}_{1}(x)\right)^{2}dx\right\},
\]
where we choose $\lambda = 3000$, $x_0 = 5$ and $\bar{\rho}_1 = \mathcal{N}(6,1)$.

The discretized loss function is 
\[
L=\frac{1}{M}\sum_{m=1}^{M}\left[\sum_{n=1}^{N}\left(\alpha_{n,m} -0.5\right)^2 \Delta t\right]+\lambda\sum_{i=1}^{I}\frac{1}{2}\left(\tilde{\rho}_{N}(x_{i})-\tilde{\bar{\rho}}_{N}(x_{i})\right)^{2}\Delta x.
\]

\begin{figure}[h]
\subfloat[compare $\tilde{\rho}_{1}$ with $\tilde{\bar{\rho}}_{1}$\label{fig:L2_compare-with-KDE}]{\includegraphics[scale=0.4]{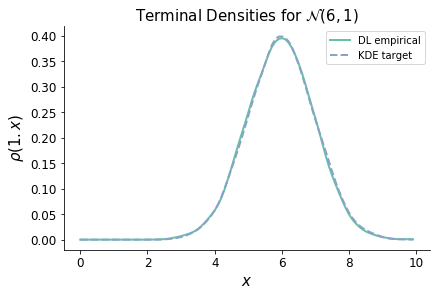}}
\subfloat[compare $\tilde{\rho}_{1}$ with $\bar{\rho}_{1}$\label{fig:L2_compare-with-true}]
{\includegraphics[scale=0.4]{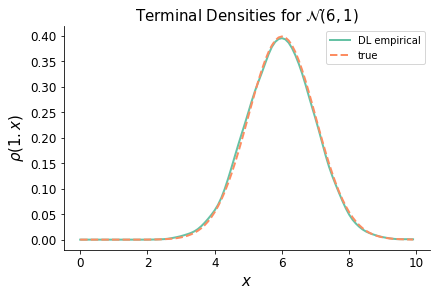}

}

\subfloat[Loss function during the training\label{fig:L2_Loss-function}]{\includegraphics[scale=0.4]{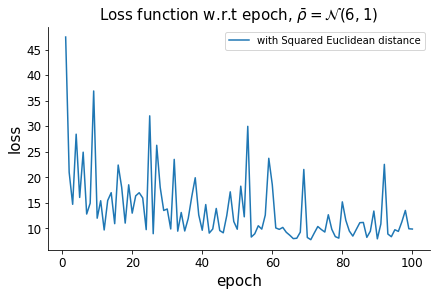}

}

{

\caption{Squared Euclidean: $\bar{\rho}=\mathcal{N}(6,1),\lambda=3000, \alpha_t \in \mathbb{R}$}

}
\end{figure}

We compare the KDE-estimated empirical terminal density $\tilde{\rho}_{1}$
with the KDE-estimated target density $\tilde{\bar{\rho}}_{1}$ in
Figure \ref{fig:L2_compare-with-KDE}; then we also compare $\tilde{\rho}_{1}$
with the true target density $\bar{\rho}_{1}$ in Figure \ref{fig:L2_compare-with-true}. We can see that the portfolio allocation learned with the deep neural network can successfully steer the density to the prescribed one. We show the change of loss function through the training process
in Figure \ref{fig:L2_Loss-function}. In training, we use a mini-batch
gradient descent algorithm, which is a mix of batch gradient descent
and stochastic gradient descent. In each step, it uses one mini-batch
to compute the gradient and update the loss function. Therefore, as
we can see in Figure \ref{fig:L2_Loss-function}, the loss function
is not monotonic, but the overall trend is decreasing.

We can see the empirical density $\tilde{\rho}_{1}$ converges to
the target $\tilde{\bar{\rho}}_{1}$ after the training. The current
result is a trade-off between the cost function in $\alpha_{t}$ and the
penalty functional over the terminal densities. In theory, we should
set $\lambda=+\infty$ in the penalty functional to make $\rho_{1}=\bar{\rho}_{1}$.
Empirically, it is sufficient to set $\lambda$ to be a large number for a similar
effect.

The distance between $\tilde{\rho}_{1}$ and $\tilde{\bar{\rho}}_{1}$
comes from the imperfect training of the neural networks. There are
two sources of error when we compare $\tilde{\rho}_{1}$ with the
true target density $\bar{\rho}_{1}$. One is the imperfectly trained
model, and the other is the KDE estimation. To reduce the first kind
of error, we can try different hyper-parameters in the neural network,
e.g., different mini-batch sizes, different structures of the network,
etc. On the other hand, we can reduce the second kind of error
by choosing more appropriate kernels and bandwidth $h$ for the KDE
method.

\begin{figure}
\subfloat[compare $\tilde{\rho}_{1}$ with $\tilde{\bar{\rho}}_{1}$]{\includegraphics[scale=0.4]{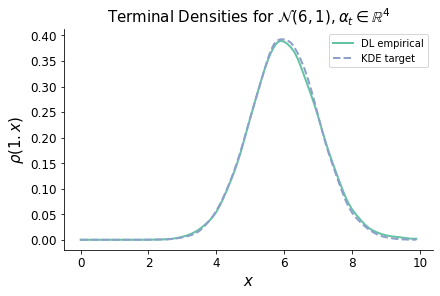}

}\subfloat[compare $\tilde{\rho}_{1}$ with $\bar{\rho}_{1}$]{\includegraphics[scale=0.4]{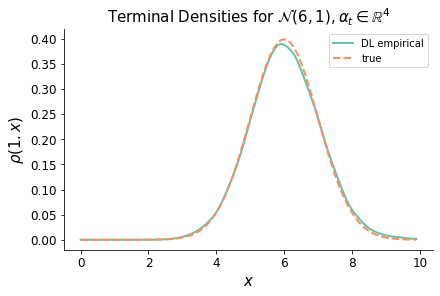}

}

{\caption{Squared Euclidean: $\bar{\rho}=\mathcal{N}(6,1),\lambda=3000,\alpha_{t}\in\mathbb{R}^{4}$\label{fig:L2_4_d}}
} 
\end{figure}

An advantage of this deep learning method is that it does not require
much extra effort for multi-asset problems. We can easily obtain
result for a 4-risky-asset case (or more) with good accuracy within computing time approximately
the same as for the 1-risky-asset case. The result for such multi-asset example
is shown in Figure \ref{fig:L2_4_d}. 

\subsubsection{Kullback-Leibler divergence}

In this second example, we use the Kullback--Leibler divergence to
measure the difference between the two distributions. 
The Kullback--Leibler (K--L) divergence (\citealt{kullback1951information}) is also known as relative entropy or information deviation.
In Machine Learning and neuroscience, the K--L divergence plays a leading role and is a widely used tool in pattern recognition (\citealt{mesaros2007singer}), multimedia classification (\citealt{moreno2004kullback}), text classification (\citealt{dhillon2003divisive}) and so on.  
In the implementation, we may face \textit{$0\log0$} or \textit{division by zero} cases in
practice; to address this, we can replace zero with an infinitesimal positive value.
In this case, the value function is defined as
\[
V(\rho_{0},\bar{\rho}_{1})=\inf_{\alpha}\left\{ \mathbb{E}\Bigl[\int_{0}^{1} (\alpha_{t}-0.5)^2 dt\Bigr]+\lambda\int_{\mathbb{R}}\rho_{1}(x)\log\left(\frac{\rho_{1}(x)}{\bar{\rho}_{1}(x)}\right)dx\right\},
\]
and the discretized form of the loss function is 
\[
L=\frac{1}{M}\sum_{m=1}^{M}\left[\sum_{n=1}^{N}\left(\alpha_{n,m}- 0.5 \right)^2\Delta t\right]+\lambda\sum_{i=1}^{I}\tilde{\rho}_{N}(x_{i})\log\left(\frac{\tilde{\rho}_{N}(x_{i})}{\tilde{\bar{\rho}}_{N}(x_{i})}\right)\Delta x.
\]

\begin{figure}
\subfloat[compare $\tilde{\rho}_{1}$ with $\tilde{\bar{\rho}}_{1}$\label{fig:KL_compare-with-KDE}]{\includegraphics[scale=0.4]{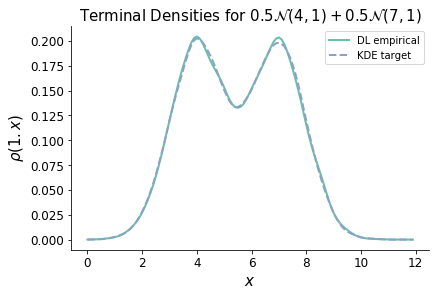}

}\subfloat[compare $\tilde{\rho}_{1}$ with $\bar{\rho}_{1}$\label{fig:KL_compare-with-true}]{\includegraphics[scale=0.4]{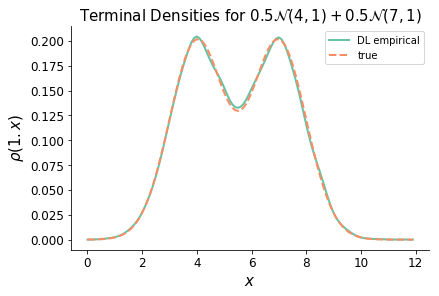}

}

\subfloat[Loss function during the training\label{fig:KL_Loss-function}]{\includegraphics[scale=0.4]{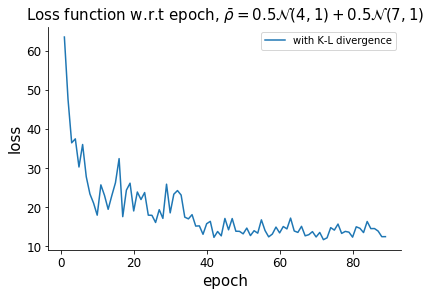}

}

\caption{K-L divergence: $\bar{\rho}_{1}=0.5\mathcal{N}(4,1)+0.5\mathcal{N}(7,1),\lambda=2000$\label{fig:K-L}}
\end{figure}

To show that our method is not restricted to Gaussian target distributions, we use a mixture of two normal distributions as the target $\bar{\rho}_{1}=0.5\mathcal{N}(4,1)+0.5\mathcal{N}(7,1)$
in this example. We provide the comparisons between
$\tilde{\rho}_{1}$ and $\tilde{\bar{\rho}}_{1}$, $\tilde{\rho}_{1}$
and $\bar{\rho}_{1}$ in figures \ref{fig:KL_compare-with-KDE} and
\ref{fig:KL_compare-with-true}, respectively. The loss function is shown afterwards, in Figure \ref{fig:KL_Loss-function}.

We want to emphasize that our method is not restricted to constant
parameters. The risky assets can have time-varying returns and covariances.
We are not illustrating the results here because the output plots
look very similar to the presented cases.

\subsubsection{2-Wasserstein Distance}

Here, we use 2-Wasserstein distance as the penalty functional. The Wasserstein distance is widely used in Machine Learning to formulate a metric for comparing clusters \citep{coen2010comparing}, and has been applied to image retrieval \citep{rubner2000earth}, contour matching \citep{grauman2004fast}, and many other problems. 
The Wasserstein distance has some advantages compared to distances such as $L^2, \chi^2$ or Hellinger. 
First of all, it can capture the underlying geometry of the space, which may be ignored by the Euclidean distance.
Secondly, when we take average of different objects -- such as distributions and images -- we can get back to a similar object with the Wasserstein distance.
Thirdly, some of the above distances are sensitive to small wiggles in the distribution, but the Wasserstein distance
is insensitive to small wiggles.

The 2-Wasserstein distance between two distributions $\rho_{1}$ and $\bar{\rho}_{1}$
is defined by
\[
W_{2}(\rho_{1},\bar{\rho}_{1})=\left(\inf_{\gamma\in\Gamma(\rho_{1},\bar{\rho}_{1})}\int_{\mathbb{R}\times\mathbb{R}}\left|x-y\right|^{2}d\gamma(x,y)\right)^{\frac{1}{2}}.
\]
For continuous one-dimensional probability distributions $\rho_{1}$
and $\bar{\rho}_{1}$ on $\mathbb{R}$, the distance has a closed
form in terms of the corresponding cumulative distribution functions
$F(x)$ and $\bar{F}(x)$ (see \citealt{ruschendorf1985wasserstein}
for the detailed proof): 
\[
W_{2}(\rho_{1},\bar{\rho}_{1})=\left(\int_{0}^{1}\left|F^{-1}(u)-\bar{F}^{-1}(u)\right|^{2}du\right)^{\frac{1}{2}}.
\]
Empirically, we can use the order statistics: 
\[
W_{2}(\rho_{1},\bar{\rho}_{1})=\left(\sum_{i=1}^{n}\left|X_{i}-Y_{i}\right|^{2}\right)^{\frac{1}{2}},
\]
where the dataset $X_{1},X_{2},...X_{n}$ is increasingly ordered
with an empirical distribution $\rho_{1}$, similarly the dataset
$Y_{1},Y_{2},...Y_{n}$ is increasingly ordered with an empirical
distribution $\bar{\rho}_{1}$. 

\begin{figure}[H]
\subfloat[compare $\tilde{\rho}_{1}$ with $\tilde{\bar{\rho}}_{1}$]{\includegraphics[scale=0.4]{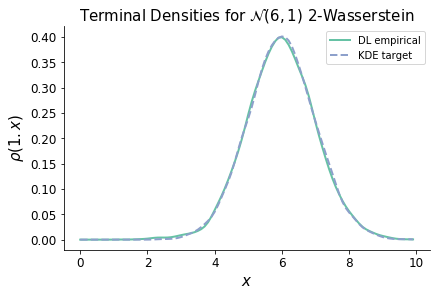}

}\subfloat[compare $\tilde{\rho}_{1}$ with $\bar{\rho}_{1}$]{\includegraphics[scale=0.4]{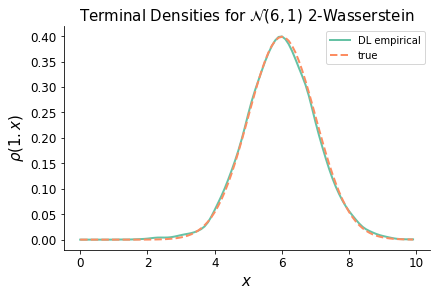}

}

\subfloat[Loss function during the training]{\includegraphics[scale=0.4]{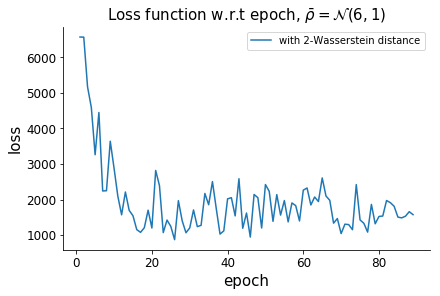}

}

\caption{2-Wasserstein distance: $\bar{\rho}_{1}=\mathcal{N}(6,1),\lambda=4000$}
\label{fig: 2_Wasserstein}
\end{figure}

Using the same cost function and network structure as before, 
the numerical results for the densities and loss function are illustrated in Figure \ref{fig: 2_Wasserstein}.

\section{Conclusion}
In this paper, we first devise a deep learning method to solve the optimal transport problem via a penalization method. In particular, we relax the classical optimal transport problem and introduce a functional to penalize the deviation between the empirical terminal density and the prescribed one. In \prettyref{sec:Application-in-Portfolio}, we apply this deep learning method to the portfolio allocation problem raised in \citet{guo2020portfolio}, where the goal is to reach a prescribed wealth distribution at the final time. We then provide numerical results for various choices of penalty functionals and target densities.

Then we investigate the dual representation and find it can be written as a saddle point problem. In this way, the optimal transport problem can be interpreted as a minimax game between the drift/diffusion and the potential function. We then solve this differential game with adversarial networks. Because this algorithm is free of spatial discretization,  it can be applied to high dimensional optimal transport problems, where we illustrate an example up to dimension $10$. These examples validate the accuracy and flexibility of the proposed deep optimal transport method.

 \bibliographystyle{chicago}
\bibliography{deep_learning_OT}

\begin{thebibliography}{}

\bibitem[\protect\citeauthoryear{Bachouch, Hur{\'e}, Langren{\'e}, and
  Pham}{Bachouch et~al.}{2021}]{bachouch2018deep}
Bachouch, A., C.~Hur{\'e}, N.~Langren{\'e}, and H.~Pham (2021).
\newblock Deep neural networks algorithms for stochastic control problems on
  finite horizon: numerical applications.
\newblock {\em Methodology and Computing in Applied Probability\/}.
\newblock To appear.

\bibitem[\protect\citeauthoryear{Chen, Pelger, and Zhu}{Chen
  et~al.}{2019}]{chen2019deep}
Chen, L., M.~Pelger, and J.~Zhu (2019).
\newblock Deep learning in asset pricing.
\newblock {\em Available at SSRN 3350138\/}.

\bibitem[\protect\citeauthoryear{Chow, Li, Osher, and Yin}{Chow
  et~al.}{2019}]{chow2019algorithm}
Chow, Y.~T., W.~Li, S.~Osher, and W.~Yin (2019).
\newblock Algorithm for {H}amilton--{J}acobi equations in density space via a
  generalized {H}opf formula.
\newblock {\em Journal of Scientific Computing\/}~{\em 80\/}(2), 1195--1239.

\bibitem[\protect\citeauthoryear{Coen, Ansari, and Fillmore}{Coen
  et~al.}{2010}]{coen2010comparing}
Coen, M.~H., M.~H. Ansari, and N.~Fillmore (2010).
\newblock Comparing clusterings in space.
\newblock In {\em Proceedings of the 27th International Conference on Machine
  Learning (ICML-10)}, pp.\  231--238.

\bibitem[\protect\citeauthoryear{Dhillon, Mallela, and Kumar}{Dhillon
  et~al.}{2003}]{dhillon2003divisive}
Dhillon, I.~S., S.~Mallela, and R.~Kumar (2003).
\newblock A divisive information-theoretic feature clustering algorithm for
  text classification.
\newblock {\em Journal of machine learning research\/}~{\em 3\/}(Mar),
  1265--1287.

\bibitem[\protect\citeauthoryear{Eckstein and Kupper}{Eckstein and
  Kupper}{2019}]{eckstein2019computation}
Eckstein, S. and M.~Kupper (2019).
\newblock Computation of optimal transport and related hedging problems via
  penalization and neural networks.
\newblock {\em Applied Mathematics \& Optimization\/}, 1--29.

\bibitem[\protect\citeauthoryear{Goodfellow, Pouget-Abadie, Mirza, Xu,
  Warde-Farley, Ozair, Courville, and Bengio}{Goodfellow
  et~al.}{2014}]{goodfellow2014generative}
Goodfellow, I., J.~Pouget-Abadie, M.~Mirza, B.~Xu, D.~Warde-Farley, S.~Ozair,
  A.~Courville, and Y.~Bengio (2014).
\newblock Generative adversarial nets.
\newblock In {\em Advances in neural information processing systems},
  Volume~27, pp.\  2672--2680.

\bibitem[\protect\citeauthoryear{Grauman and Darrell}{Grauman and
  Darrell}{2004}]{grauman2004fast}
Grauman, K. and T.~Darrell (2004).
\newblock Fast contour matching using approximate {E}arth {M}over's {D}istance.
\newblock In {\em Proceedings of the 2004 IEEE Computer Society Conference on
  Computer Vision and Pattern Recognition, 2004. CVPR 2004.}, Volume~1, pp.\
  I--I. IEEE.

\bibitem[\protect\citeauthoryear{Guo, Langren{\'e}, Loeper, and Ning}{Guo
  et~al.}{2019}]{guo2019robust}
Guo, I., N.~Langren{\'e}, G.~Loeper, and W.~Ning (2019).
\newblock Robust utility maximization under model uncertainty via a
  penalization approach.
\newblock {\em arXiv preprint arXiv:1907.13345\/}.

\bibitem[\protect\citeauthoryear{Guo, Langren{\'e}, Loeper, and Ning}{Guo
  et~al.}{2020}]{guo2020portfolio}
Guo, I., N.~Langren{\'e}, G.~Loeper, and W.~Ning (2020).
\newblock Portfolio optimization with a prescribed terminal wealth
  distribution.
\newblock {\em arXiv preprint arXiv:2009.12823\/}.

\bibitem[\protect\citeauthoryear{Haber and Horesh}{Haber and
  Horesh}{2015}]{haber2015multilevel}
Haber, E. and R.~Horesh (2015).
\newblock A multilevel method for the solution of time dependent optimal
  transport.
\newblock {\em Numerical Mathematics: Theory, Methods and Applications\/}~{\em
  8\/}(1), 97--111.

\bibitem[\protect\citeauthoryear{Han and E}{Han and E}{2016}]{han2016deep}
Han, J. and W.~E (2016).
\newblock Deep learning approximation for stochastic control problems.
\newblock In {\em NIPS 2016, Deep Reinforcement Learning Workshop}.

\bibitem[\protect\citeauthoryear{Han and Jentzen}{Han and
  Jentzen}{2020}]{han2020algorithms}
Han, J. and A.~Jentzen (2020).
\newblock Algorithms for solving high dimensional {PDE}s: from nonlinear
  {M}onte {C}arlo to {M}achine {L}earning.
\newblock {\em arXiv preprint arXiv:2008.13333\/}.

\bibitem[\protect\citeauthoryear{Henry-Labord\`ere}{Henry-Labord\`ere}{2019}]{henry2019martingale}
Henry-Labord\`ere, P. (2019).
\newblock ({M}artingale) optimal transport and anomaly detection with neural
  networks: a primal-dual algorithm.
\newblock {\em Available at SSRN 3370910\/}.

\bibitem[\protect\citeauthoryear{Hur{\'e}, Pham, Bachouch, and
  Langren{\'e}}{Hur{\'e} et~al.}{2021}]{hure2018deep}
Hur{\'e}, C., H.~Pham, A.~Bachouch, and N.~Langren{\'e} (2021).
\newblock Deep neural networks algorithms for stochastic control problems on
  finite horizon: convergence analysis.
\newblock {\em SIAM Journal on Numerical Analysis\/}~{\em 59\/}(1), 525--557.

\bibitem[\protect\citeauthoryear{Hur{\'e}, Pham, and Warin}{Hur{\'e}
  et~al.}{2019}]{hure2019some}
Hur{\'e}, C., H.~Pham, and X.~Warin (2019).
\newblock Some machine learning schemes for high-dimensional nonlinear {PDE}s.
\newblock {\em arXiv preprint arXiv:1902.01599\/}.

\bibitem[\protect\citeauthoryear{Kantorovich}{Kantorovich}{1942}]{kantorovich1942translocation}
Kantorovich, L.~V. (1942).
\newblock On the translocation of masses.
\newblock In {\em Dokl. Akad. Nauk. USSR (NS)}, Volume~37, pp.\  199--201.

\bibitem[\protect\citeauthoryear{Kullback and Leibler}{Kullback and
  Leibler}{1951}]{kullback1951information}
Kullback, S. and R.~A. Leibler (1951).
\newblock On information and sufficiency.
\newblock {\em The Annals of Mathematical Statistics\/}~{\em 22\/}(1), 79--86.

\bibitem[\protect\citeauthoryear{Li, Ryu, Osher, Yin, and Gangbo}{Li
  et~al.}{2018}]{li2018parallel}
Li, W., E.~K. Ryu, S.~Osher, W.~Yin, and W.~Gangbo (2018).
\newblock A parallel method for {E}arth {M}over's {D}istance.
\newblock {\em Journal of Scientific Computing\/}~{\em 75\/}(1), 182--197.

\bibitem[\protect\citeauthoryear{Liang and Srikant}{Liang and
  Srikant}{2016}]{liang2016deep}
Liang, S. and R.~Srikant (2016).
\newblock Why deep neural networks for function approximation?
\newblock {\em arXiv preprint arXiv:1610.04161\/}.

\bibitem[\protect\citeauthoryear{Loeper}{Loeper}{2006}]{loeper2006reconstruction}
Loeper, G. (2006).
\newblock The reconstruction problem for the {E}uler-{P}oisson system in
  cosmology.
\newblock {\em Archive for rational mechanics and analysis\/}~{\em 179\/}(2),
  153--216.

\bibitem[\protect\citeauthoryear{McCulloch and Pitts}{McCulloch and
  Pitts}{1943}]{mcculloch1943logical}
McCulloch, W.~S. and W.~Pitts (1943).
\newblock A logical calculus of the ideas immanent in nervous activity.
\newblock {\em The Bulletin of Mathematical Biophysics\/}~{\em 5\/}(4),
  115--133.

\bibitem[\protect\citeauthoryear{Mesaros, Virtanen, and Klapuri}{Mesaros
  et~al.}{2007}]{mesaros2007singer}
Mesaros, A., T.~Virtanen, and A.~Klapuri (2007).
\newblock Singer identification in polyphonic music using vocal separation and
  pattern recognition methods.
\newblock In {\em ISMIR}, pp.\  375--378.

\bibitem[\protect\citeauthoryear{Monge}{Monge}{1781}]{monge1781memoire}
Monge, G. (1781).
\newblock M{\'e}moire sur la th{\'e}orie des d{\'e}blais et des remblais.
\newblock {\em Histoire de l'Acad{\'e}mie Royale des Sciences de Paris\/}.

\bibitem[\protect\citeauthoryear{Moreno, Ho, and Vasconcelos}{Moreno
  et~al.}{2004}]{moreno2004kullback}
Moreno, P.~J., P.~P. Ho, and N.~Vasconcelos (2004).
\newblock A {K}ullback-{L}eibler divergence based kernel for {SVM}
  classification in multimedia applications.
\newblock In {\em Advances in neural information processing systems}, pp.\
  1385--1392.

\bibitem[\protect\citeauthoryear{Peyr{\'e} and Cuturi}{Peyr{\'e} and
  Cuturi}{2019}]{peyre2019computational}
Peyr{\'e}, G. and M.~Cuturi (2019).
\newblock Computational optimal transport: with applications to data science.
\newblock {\em Foundations and Trends{\textregistered} in Machine
  Learning\/}~{\em 11\/}(5-6), 355--607.

\bibitem[\protect\citeauthoryear{Rachev and R{\"u}schendorf}{Rachev and
  R{\"u}schendorf}{1998}]{rachev1998mass}
Rachev, S.~T. and L.~R{\"u}schendorf (1998).
\newblock {\em Mass {T}ransportation {P}roblems: {V}olume {I}: {T}heory},
  Volume~1.
\newblock Springer Science \& Business Media.

\bibitem[\protect\citeauthoryear{Rubner, Tomasi, and Guibas}{Rubner
  et~al.}{2000}]{rubner2000earth}
Rubner, Y., C.~Tomasi, and L.~J. Guibas (2000).
\newblock The {E}arth {M}over's {D}istance as a metric for image retrieval.
\newblock {\em International journal of computer vision\/}~{\em 40\/}(2),
  99--121.

\bibitem[\protect\citeauthoryear{R{\"u}schendorf}{R{\"u}schendorf}{1985}]{ruschendorf1985wasserstein}
R{\"u}schendorf, L. (1985).
\newblock The {W}asserstein distance and approximation theorems.
\newblock {\em Probability Theory and Related Fields\/}~{\em 70\/}(1),
  117--129.

\bibitem[\protect\citeauthoryear{Ruthotto, Osher, Li, Nurbekyan, and
  Fung}{Ruthotto et~al.}{2020}]{ruthotto2020machine}
Ruthotto, L., S.~J. Osher, W.~Li, L.~Nurbekyan, and S.~W. Fung (2020).
\newblock A machine learning framework for solving high-dimensional mean field
  game and mean field control problems.
\newblock {\em Proceedings of the National Academy of Sciences\/}~{\em
  117\/}(17), 9183--9193.

\bibitem[\protect\citeauthoryear{Santambrogio}{Santambrogio}{2015}]{santambrogio2015optimal}
Santambrogio, F. (2015).
\newblock Optimal transport for applied mathematicians.
\newblock {\em Birk{\"a}user, NY\/}~{\em 55\/}(58-63), 94.

\bibitem[\protect\citeauthoryear{Sirignano and Spiliopoulos}{Sirignano and
  Spiliopoulos}{2018}]{sirignano2018dgm}
Sirignano, J. and K.~Spiliopoulos (2018).
\newblock {DGM}: A deep learning algorithm for solving partial differential
  equations.
\newblock {\em Journal of Computational Physics\/}~{\em 375}, 1339--1364.

\bibitem[\protect\citeauthoryear{Tan and Touzi}{Tan and
  Touzi}{2013}]{tan2013optimal}
Tan, X. and N.~Touzi (2013).
\newblock Optimal transportation under controlled stochastic dynamics.
\newblock {\em The Annals of Probability\/}~{\em 41\/}(5), 3201--3240.

\bibitem[\protect\citeauthoryear{Villani}{Villani}{2003}]{villani2003topics}
Villani, C. (2003).
\newblock {\em Topics in optimal transportation}, Volume~58 of {\em Graduate
  Studies in Mathematics}.
\newblock American Mathematical Society.

\bibitem[\protect\citeauthoryear{Villani}{Villani}{2008}]{villani2008optimal}
Villani, C. (2008).
\newblock {\em Optimal transport: old and new}, Volume 338 of {\em Grundlehren
  der mathematischen Wissenschaften}.
\newblock Springer Science \& Business Media.

\bibitem[\protect\citeauthoryear{Weinan, Han, and Jentzen}{Weinan
  et~al.}{2017}]{weinan2017deep}
Weinan, E., J.~Han, and A.~Jentzen (2017).
\newblock Deep learning-based numerical methods for high-dimensional parabolic
  partial differential equations and backward stochastic differential
  equations.
\newblock {\em Communications in Mathematics and Statistics\/}~{\em 5\/}(4),
  349--380.

\bibitem[\protect\citeauthoryear{Wiese, Knobloch, Korn, and Kretschmer}{Wiese
  et~al.}{2020}]{wiese2020quant}
Wiese, M., R.~Knobloch, R.~Korn, and P.~Kretschmer (2020).
\newblock Quant {GAN}s: deep generation of financial time series.
\newblock {\em Quantitative Finance\/}, 1--22.

\bibitem[\protect\citeauthoryear{Zhang, Zhong, and Ma}{Zhang
  et~al.}{2020}]{zhang2020review}
Zhang, J., W.~Zhong, and P.~Ma (2020).
\newblock A review on modern computational optimal transport methods with
  applications in biomedical research.
\newblock {\em arXiv preprint arXiv:2008.02995\/}.

\bibitem[\protect\citeauthoryear{Zhang, Zhong, Dong, Wang, and Wang}{Zhang
  et~al.}{2019}]{zhang2019stock}
Zhang, K., G.~Zhong, J.~Dong, S.~Wang, and Y.~Wang (2019).
\newblock Stock market prediction based on generative adversarial network.
\newblock {\em Procedia computer science\/}~{\em 147}, 400--406.

\end{thebibliography}

\end{document}